\input amstex
\documentstyle{amsppt}

\input epsf

\pagewidth{30pc}
\pageheight{48.6pc}

\catcode`\@=11
\let\logo@=\relax
\catcode`\@=13

\loadbold

\define\A{\text{\rm A}}
\redefine\B{\text{\rm B}}
\define\C{\text{\rm C}}
\redefine\D{\text{\rm D}}
\define\E{\text{\rm E}}
\define\F{\text{\rm F}}
\define\R{\text{\rm R}}
\define\V{\text{\rm V}}

\define\etal{\text{\it et al\. }}
\define\mc{\Cal}
\define\mb{\bold}
\define\mbb{\Bbb}
\define\mf{\frak}
\define\ellom{\ell_{\omega}}
\define\wh{\widehat}
\define\wt{\widetilde}
\define\kpa{\kappa}
\define\von{\varepsilon}
\define\tta{\theta}
\define\lba{\lambda}
\define\raro{\rightarrow}
\define\laro{\leftarrow}
\define\sumin{\sum_{i=1}^{n}}
\define\setR{\Bbb R}
\define\setN{\Bbb N}
\define\setRp{\Bbb R^{+}}
\define\GV{\Gamma\V}
\define\RV{\R\V}
\define\PV{\Pi\V}
\define\eqdef{\triangleq}
\define\trans{{}^{t}}
\define\Xn{\frak X_n}
\define\Prob{\Bbb P}
\define\Pm{\Prob_{\mu}^{n}}
\define\Pfm{\Prob_{f,\mu}^{n}}
\define\Efm{\Bbb E_{\,f, \mu}^{\,n}}
\define\Em{\Bbb E_{\mu}^{n}}

\define\mbt#1{\wt{\mb #1}}
\define\ind#1{{\bold 1}_{#1}}
\define\prodsca#1#2{\langle{#1},{#2}\rangle}
\define\prodscahk#1#2{\prodsca{#1}{#2}_{h,K}}
\define\norminfty#1{\|#1\|_{\infty}}
\define\norm#1{\|#1\|}
\define\normhk#1{\norm{#1}_{h,K}}
\define\ppint#1{\lfloor#1\rfloor}

\define\supp{\operatorname{Supp}}
\define\argmin{\operatornamewithlimits{argmin}}

\def\zs#1{_{\lower2pt\hbox{$\scriptstyle#1$}}}

\leftheadtext\nofrills{\eightpoint \it
S.~Ga{\"\i}ffas}
\rightheadtext\nofrills{\eightpoint \it
Convergence Rates with a Degenerate Design}

\line{\eightpoint Volume ~14, No.~1~(2005), pp\.
\hfil Allerton Press, Inc\.}
\smallskip
\hrule\smallskip\hrule
\medskip
\line {$\Bbb {M\,\,A\,\,T\,\,H\,\,E\,\,M\,\,A\,\,T\,\,I\,\,C\,\,A\,\,L}$\hfil
$\Bbb{M\,\,E\,\,T\,\,H\,\,O\,\,D\,\,S}$\hfil
$\Bbb{O\,\,F}$\hfil
$\Bbb{S\,\,T\,\,A\,\,T\,\,I\,\,S\,\,T\,\,I\,\,C\,\,S}$}
\medskip
\hrule\smallskip\hrule

\medskip\medskip\medskip\medskip\medskip
\medskip\medskip\medskip\medskip\medskip

\centerline{\bf CONVERGENCE RATES FOR POINTWISE CURVE}
\centerline{\bf ESTIMATION WITH A DEGENERATE DESIGN}
\medskip\medskip
\centerline{\tensmc S.~Ga{\"\i}ffas%
\footnote""{\copyright 2005 by
Allerton Press, Inc\. Authorization to photocopy individual
items for internal or personal use, or the internal or personal
use of specific clients, is granted by Allerton Press, Inc\. for
libraries and other users registered with the Copyright
Clearance Center (CCC) Transactional Reporting Service, provided
that the base fee of \$50.00 per copy is paid directly to CCC,
222 Rosewood Drive, Danvers, MA 01923.}}

{\eightpoint

\medskip\medskip
\centerline{Labor\. Probab\. et Mod{\`e}les Al{\'e}atoires,
U.M.R. CNRS 7599 and Univ\. Paris 7}
\centerline{175 rue du Chevaleret, 75013 Paris}
\centerline{E-mail: gaiffas\@math.jussieu.fr}
\medskip\medskip

}

\hrule\medskip
\par{\eightpoint\sl
The nonparametric regression with a random design model is
considered. We want to recover the regression function at a
point $x_0$ where the design density is vanishing or exploding.
Depending on assumptions on local regularity of the regression
function and on the local behaviour of the design,
we find several minimax rates.
These rates lie in a wide range, from slow $\ell(n)$ rates,
where $\ell$ is slowly varying
{\rm(}for instance $(\log n)^{-1}${\rm)},
to fast $n^{-1/2} \ell(n)$ rates. If the continuity modulus of
the regression function at $x_0$ can be bounded from above by an
$s$-regularly varying function, and if the design density is
$\beta$-regularly varying, we prove that the minimax convergence
rate at $x_0$ is $n^{-s/(1+2s+\beta)}\ell(n)$.
\medskip

\par
Key words:
degenerate design, minimax, nonparametric regression,
random design
\smallskip

\par
2000 Mathematics Subject Classification: 62G05, 62G08.}
\medskip
\hrule
\medskip\medskip\medskip\medskip\medskip

{\bf 1. Introduction}

\qquad{\tensmc1.1. The model}.
Suppose that we have $n$ independent and identically distributed
observations $(X_i, Y_i) \in \setR \times \setR$ from the
regression model
$$
Y_i = f(X_i) + \xi_i,
\tag 1.1
$$
where $f \: \setR \raro \setR$, the variables $(\xi_i)$ are
centered Gaussian of variance $\sigma^2$ and independent of $X_1,
\ldots, X_n$ (the design), and the $X_i$ are distributed with
density~$\mu$. We want to recover $f$ at a chosen~$x_0$.

For instance, if we take the variables $(X_i)$ distributed with
density
$$
\mu(x) = \frac{\beta+1}{x_0^{\beta+1} + (1-x_0)^{\beta+1}}
|x-x_0|^{\beta} \ind{[0,1]}(x),
$$
for $x_0 \in [0,1]$ and $\beta > -1$, then clearly when $\beta > 0$
this density models a lack of information at $x_0$ and conversely an
exploding amount of information if $-1<\beta<0$.  We want to
understand the influence of the parameter $\beta$ on the amount of
information at~$x_0$ in the minimax setup.  \medskip

\qquad{\tensmc1.2. Motivations}.  The pointwise estimation of the
regression function is a well-known problem, which has been
intensively studied by many authors. The first authors who computed
the minimax rate over a nonparametric class of H{\"o}lderian functions
were Ibragimov and Hasminski~(1981) and Stone~(1977).  Over the class
of H{\"o}lder functions with smoothness~$s$, the local polynomial
estimator converges with the rate $n^{-s/(1+2s)}$ (see~ Stone~(1977))
and this rate is optimal in the minimax sense. Many authors worked on
related problems: see, for instance, Korostelev and Tsybakov~(1993),
Nemirovski~(2000), Tsybakov~(2003).

Nevertheless, these results require the design density to be
non-vanishing and finite at the estimation point. This assumption
roughly means that the information is spatially {\it homogeneous}.
The next logical step is to look for the minimax risk at a point where
the design density $\mu$ is vanishing or exploding. To achieve such a
result, it seems natural to consider several types of design
density behaviour at $x_0$ and to compute the corresponding
minimax rates.  Such results would improve the statistical
description of models (here in the minimax setup) with very
inhomogeneous information.

When $f$ has a H{\"o}lder type smoothness of order~$2$ and if $\mu(x)
\sim x ^{\beta}$ near $0$, where $\beta > 0$, Hall \etal (1997) show
that a local linear procedure converges with the rate $n^{-4 / (5 +
  \beta)}$ when estimating $f$ at~$0$.  This rate is also proved to be
optimal. In a more general setup for the design and if the regression
function is Lipschitz, Guerre~(1999) extends the result of Hall \etal
for $\beta>-1$.  Here, we intend to develop the regression function
estimation for degenerate designs in a systematic way.
\medskip

\qquad{\tensmc1.3. Organization of the paper}.
In Section~2 we present two theorems giving the pointwise
minimax convergence rate in the model (1.1) for different design
behaviours (Theorems~1 and~2).
In Section~3 we construct an estimator and in Section~4
give upper bounds for this estimator
(Propositions~4 and~5).
In Section~5 we discuss some technical points. The proofs are
delayed until Section~6 and well-known facts about the regular
and $\Gamma$-variation are given in the Appendix.
\medskip

{\bf 2. Main Results}

All along this study we are in the minimax setup. We define the
pointwise minimax risk over a class~$\Sigma$ by
$$
\mc R_n(\Sigma, \mu) \eqdef
\Bigl(\inf_{T_n} \sup_{f \in \Sigma}
\Efm \{ | T_n(x_0) - f(x_0) |^p \} \Bigr)^{1/p},
\tag 2.1
$$
where $\inf_{T_n}$ is taken over all estimators $T_n$ based on the
observations (1.1), with $x_0$ being the estimation point and
$p>0$. The expectation $\Efm$ in (2.1) is taken with respect to
the joint probability distribution $\Pfm$ of the pairs
$(X_i,Y_i)_{i=1,\ldots,n}$.
\medskip

\qquad{\tensmc2.1. Regular variation}.
The definition of regular variation and the main properties are
due to Karamata (1930). The main references on regular variation
are Bingham \etal (1989), Geluk and de Haan~(1987),
Resnick (1987), and Senata~(1976).
\medskip

{\bf Definition 1} (Regular variation).
A continuous function $\nu \: \setRp \rightarrow \setRp$ is
regularly varying at~$0$ if there is a real number
$\beta \in \setR$ such that:
$$
\forall y > 0, \quad
\lim_{h \rightarrow 0^+} \nu(yh) / \nu(h) = y^{\beta}.
\tag 2.2
$$
We denote by $\RV(\beta)$ the set of all the functions
satisfying (2.2). A function in $\RV(0)$ is
{\it slowly varying}.
\medskip

{\bf Remark.}
Roughly, a regularly varying function behaves as a power
function times a slower term. Typical examples of such functions
are $x^{\beta}$, $x^{\beta}(\log(1/x))^{\gamma}$
for $\gamma \in\setR$, and more generally any power function
times a $\log$ or a composition of $\log$-functions to some
power. For other examples, see the references cited above.
\medskip

\qquad{\tensmc2.2. The functions class}
\smallskip

{\bf Definition 2.}
If $\delta > 0$ and $\omega \in \RV(s)$ with $s > 0$ we define
the  class $\mc F_{\delta}(x_0, \omega)$ of functions
$f \: [0,1]\raro\setR$ such that
$$
\forall h \leq \delta, \quad \inf_{P \in \mc P_{k}}
\sup_{|x-x_0| \leq h} |f(x) - P(x-x_0)| \leq \omega(h),
$$
where $k = \ppint{s}$ (the largest integer smaller than $s$) and
$\mc P_k$ is the set of all the real polynomials with degree
$k$. We define $\ellom(h) \eqdef \omega(h) h^{-s}$, the slow
variation term of~$\omega$. If $\alpha > 0$, we define
$$
\mc U(\alpha) \eqdef \bigl\{ f \: [0,1] \raro \setR
\text{ such that } \norminfty{f} \leq \alpha \bigr\}.
$$
Finally, we define
$$
\Sigma_{\delta, \alpha}(x_0, \omega)
\eqdef \mc F_{\delta}(x_0,\omega) \cap \mc U(\alpha).
$$
\medskip

{\bf Remark.}
If we take $\omega(h) = r h^s$ for some $r > 0$, then we get
the classical H{\"o}lder regularity with radius~$r$. In this
sense, the class $\mc F_{\delta}(x_0, \omega)$ is a slight
generalization of the H{\"o}lder regularity.
\medskip

{\bf Assumption M.}
In what follows, we assume that there exists a neighbourhood
$W$ of $x_0$ and a continuous function $\nu\:\setRp\raro\setRp$
such that:
$$
\forall x \in W, \quad \mu(x) = \nu(|x - x_0|).
\tag 2.3
$$
\medskip

This assumption roughly means that close to $x_0$ there are as
many observations on the left of $x_0$ as on the right. All
the following results can be extended easily to
the non-symmetric case, see Section~5.1.
\medskip

\qquad{\tensmc2.3. Regularly varying design density}.
Theorem~1 gives the minimax rate over the class $\Sigma$
(see Definition~2) for the estimation problem of $f$ at $x_0$
when the design is regularly varying at this point.

We denote by $\mc R(x_0,\beta)$ the set of all the densities
$\mu$ such that (2.3) holds with $\nu \in \RV(\beta)$ for a
fixed neighbourhood~$W$.
\medskip

{\bf Theorem 1.}
{\it If
\roster
\item"$\bullet$"
$(s,\beta) \in (0,+\infty) \times (-1, +\infty)$ or
$(s,\beta) \in (0,1] \times \{-1\}$,
\item"$\bullet$"
$\Sigma = \Sigma_{h_n, \alpha_n}(x_0, \omega)$ with
$\omega \in \RV(s)$, $\alpha_n = O(n^{\gamma})$
for some $\gamma > 0$ and~$h_n$ given by {\rm(2.5)},
\item"$\bullet$"
$\mu \in \mc R(x_0,\beta)$,
\endroster
then we have
$$
\mc R_n(\Sigma,\mu) \asymp \sigma^{2s/(1 + 2s + \beta)}
n^{-s/(1 + 2s + \beta)} \ell_{\omega,\nu}(n^{-1})\qquad
\text{as} \quad n \raro +\infty,
\tag 2.4
$$
where $\ell_{\omega,\nu}$ is slowly varying and where $\asymp$
stands for the equality in order, up to constants depending on
$s$, $\beta$ and $p$ {\rm(}see {\rm(2.1))}
but not on $\sigma$.  Moreover, the minimax rate is equal to
$\omega(h_n)$, where $h_n$ is the smallest solution to
$$
\omega(h) = \frac{\sigma}{\sqrt{2 n \int_0^h \nu(t)\, dt}}.
\tag 2.5
$$}

{\bf Example.}
The simplest example is the non-degenerate design case
($0 < \mu(x_0) <+\infty$) with the class $\Sigma$ equal to a
H\"older ball ($\omega(h) = r h^s$, see Definition~2).
This is the common case found in the literature.
In particular,
in this case, the design is slowly varying ($\beta = 0$ with the
slow term constant and equal to $\lim_{x \raro x_0} \mu(x)$).
Solving (2.5) leads to the classical minimax rate
$$
\sigma^{2s/(1+2s)} r^{1/(1+2s)} n^{-s/(1+2s)}.
$$

{\bf Example.}
Let $\beta > -1$. We consider $\nu$ such that
$\int_0^h \nu(t)\,dt = h^{\beta+1}(\log(1/h))^{\alpha}$ and
$\omega(h) = r h^s (\log(1/h))^{\gamma}$, where $\alpha, \gamma$
are any real numbers.
In this case, we find that the minimax rate
(see Section~6.5 for details) is
$$
\sigma^{2s/(1+2s+\beta)} r^{(\beta+1)/(1+2s+\beta)}
\big(n(\log n)^{\alpha - \gamma(1+\beta)/s}\big)^{-s/(1+2s+\beta)}.
$$

We note that this rate has the form given by Theorem~1
with the slow term $\ell_{\omega,\nu}(h) =
(\log(1/h))^{(\gamma(\beta+1) -s \alpha)/(1+2s+\beta)}$.
When $\gamma(1+\beta) - s\alpha = 0$, there is no slow term in
the minimax rate, although there are slow terms in $\nu$ and
$\omega$. Again, if $\beta = 0$ and $\gamma = s\alpha$, we get
the minimax rate of the first example, although the terms $\nu$
and $\omega$ do not have the classical forms.
\medskip

{\bf Example.}
Let $\beta = -1$, $\alpha>1$, and $\nu(h) = h^{-1}
(\log(1/h))^{-\alpha}$. Let $\omega$ be the same as in the
previous example with $0 < s \leq 1$. Then the minimax
convergence rate is
$$
\sigma n^{-1/2} (\log n)^{(\alpha-1)/2}.
$$

This rate is almost the parametric estimation rate, up to the
slow $\log$ factor. This result is natural since the design is
very ``exploding": we have a lot of information at $x_0$,
thus we can estimate $f(x_0)$ very fast. Also, we note that the
regularity parameters of the regression function ($r$, $s$, and
$\gamma$) have (asymptotically) disappeared from the minimax
rate.
\medskip

\qquad{\tensmc2.4. $\Gamma$-varying design density}.
The regular variation framework includes any design density
behaving close to the estimation point as a polynomial times a
slow term. It does not include, for instance, a design with a
behaviour similar to $\exp(-1/|x-x_0|)$ and defined as~$0$
at~$x_0$, since this function goes to $0$ at $x_0$ faster
than any power function.

Such a local behaviour can model the situation where we have
very little information.
This example naturally leads us to the framework of
$\Gamma$-variation. In fact, such a function belongs to the
following class introduced by de Haan (1970).
\medskip

{\bf Definition 3} ($\Gamma$-variation).
A non-decreasing continuous function
$\nu \: \setR^+ \raro  \setR^+$ is $\Gamma$-varying if there
exists a continuous function
$\rho \: \setR^+ \raro \setR^+$ such that
$$
\forall y \in \setR, \quad
\lim_{h \raro 0^+} \nu(h + y \rho(h)) / \nu(h) = \exp(y).
\tag 2.6
  $$
We denote by $\GV(\rho)$ the class of all such functions. The
function $\rho$ is called the {\it auxiliary\/} function of~$\nu$.
\medskip

{\bf Remark.}
A function behaving like $\exp(-1/|x-x_0|)$ close to $x_0$
satisfies Assumption~M with $\nu(h) = \exp(-1/h)$, where
$\nu\in \GV(\rho)$ with $\rho(h) = h^2$.
\medskip

{\bf Theorem 2.}
{\it If
\roster
\item"$\bullet$"
$\Sigma = \Sigma_{h_n,\alpha_n}(x_0, \omega)$, where
$\omega \in \RV(s)$ with $0 < s \leq 1$, $h_n$ is given by
{\rm(2.5)} and $\alpha_n = O(r_n^{-\gamma})$ for some
$\gamma > 0$ with $r_n \eqdef \omega(h_n)$,
\item"$\bullet$"
$\mu$ satisfies Assumption~{\rm M} with $\nu\in\GV(\rho)$,
\endroster
then
$$
\mc R_n(\Sigma,\mu) \asymp \ell_{\omega,\nu}(n^{-1})\qquad
\text{as} \quad n \raro +\infty,
$$
where $\ell_{\omega,\nu}$ is slowly varying. Moreover, as in
Theorem~\rom1, the minimax rate is equal to $\omega(h_n)$,
where $h_n$ is the smallest solution to~{\rm(2.5)}.}
\medskip

{\bf Example.}
Let $\mu$ satisfy Assumption~M  with
$\nu(h)=\exp(-1/h^{\alpha})$ for $\alpha > 0$ and
$\omega(h) = r h^s$ for $0 < s \leq 1$.
It is an easy computation to see that $\nu$ belongs to
the class $\GV(\rho)$ for the auxiliary function
$\rho(h) =  \alpha^{-1} h^{\alpha+1}$. In this case, we find
that the minimax rate (see Section~6.5 for details) is
$$
r(\log n)^{-s/\alpha}.
$$
As shown by Theorem~2, we find a very slow minimax rate in this
example. We note that the parameters $s$ and $\alpha$ are on the
same scale.
\medskip

{\bf 3. Local Polynomial Estimation}

\qquad{\tensmc3.1. Introduction}.
For the proof of the upper bound in Theorem~1 we use
a local polynomial estimator. The local polynomial estimator is
well-known and has been intensively studied (see Stone~(1977),
Fan and Gijbels~(1996), Spokoiny~(1998),
Tsybakov~(2003), among many others). If $f$ is a smooth
function at $x_0$, then it is close to its Taylor polynomial. A
function $f\in C^k(x_0)$ (the space of $k$ times differentiable
functions at $x_0$ with a continuous $k$-th derivative) is such
that for any $x$ close to~$x_0$
$$
f(x) \approx f(x_0) + f^{'}(x_0)(x-x_0) + \ldots
+  \frac{f^{(k)}(x_0)}{k!}(x-x_0)^k.
\tag 3.1
$$
Let $h>0$ (the \text{\it bandwidth\/}) and $k \in \setN$. We
define $\phi_{j,h}(x) \eqdef \bigl(\frac{x-x_0}{h}\bigr)^j$
and the space
$$
V_{k,h} \eqdef \text{Span}\{ (\phi_{j,h})_{j=0,\ldots,k}\}.
$$
For a fixed non-negative function $K$ (the {\it kernel\/}) we
define the weighted pseudo-scalar product
$$
\prodscahk{f}{g} \eqdef \sumin f(X_i) g(X_i)
K\Bigl(\frac{X_i-x_0}{h}\Bigr),
\tag 3.2
$$
and the corresponding pseudo-norm
$\normhk{\cdot} \eqdef \sqrt{\prodscahk{\cdot}{\cdot}}$\;\
($K \geq 0$). In view of (3.1) it is natural to consider the
estimator defined as the closest polynomial of degree~$k$
to the observations $(Y_i)$ in the least square sense, that is:
$$
\wh f_{h} = \argmin_{ g \in V_{k,h}} \normhk{g - Y}^2.
\tag 3.3
$$
Then $\wh f_{h}(x_0)$ is the {\it local polynomial estimator\/}
of $f$ at $x_0$. A necessary condition for $\wh f_h$ to be the
minimizer of (3.3) is that it solves the linear problem:
$$
\text{find $\wh f \in V_{k,h}$  such that
$\forall\, \phi \in V_{k,h}$}, \quad
\prodscahk{\wh f}{\phi} = \prodscahk{Y}{\phi}.
\tag 3.4
$$
The estimator $\wh f_{h}$ is then given by
$$
\wh f_{h} = P_{\wh \theta_h},
\tag 3.5
$$
where
$$
P_{\theta} = \tta_0 \phi_{0,h} + \tta_1 \phi_{1,h} + \ldots
+  \tta_k \phi_{k,h},
\tag 3.6
$$
with $\wh \theta_h$ the solution, whenever it makes sense, of the
linear system
$$
\mb X_h^K \theta = \mb Y_h^K,
\tag 3.7
$$
where $\mb X_h^K$ is the symmetric matrix with entries
$$
(\mb X_h^K)_{j,l} = \prodscahk{\phi_{j,h}}{\phi_{l,h}},
\quad 0 \leq j,l \leq k,
\tag 3.8
$$
and $\mb Y_h^{K}$ is the vector defined by
$$
\mb Y_h^{K} = (\prodscahk{Y}{\phi_{j,h}} ; 0 \leq j \leq k ).
$$
We assume that the kernel $K$ satisfies the following assumptions:

{\bf Assumption K.}
Let $K$ be the rectangular kernel
$K^R(x) = \frac{1}{2} \ind{|x|\leq 1}$
or a non-negative function such that:
\roster
\item"$\bullet$"
$\supp{K}\subset [-1,1] $,
\item"$\bullet$"
$K$ is symmetric,
\item"$\bullet$"
$K_{\infty} \eqdef \sup_x K(x) \leq 1$,
\item"$\bullet$"
there is some $\rho>0$ and $\kappa > 0$ such that
$\forall x,y$, $|K(x)-K(y)| \leq \rho |x-y|^{\kappa}$.
\endroster
\medskip

Assumption~K is satisfied by all the classical kernels used in
nonparametric curve smoothing. Let us define
$$
N_{n,h} = \#\{X_i \text{ such that } X_i \in [x_0-h, x_0+h]\},
\tag 3.9
$$
the number of observations in the interval $[x_0-h, x_0+h]$, and
the random matrix
$$
\mc X_h^K \eqdef N_{n,h}^{-1} \mb X_h^K.
$$
Denote
$
\mf X_n \eqdef \sigma(X_1,\ldots,X_n)$
the $\sigma$-algebra generated by the design. Note that
$\mc X_h^K$ is measurable with respect to $\mf X_n$. The matrix
$\mc X_h^K$ is a ``renormalization" of $\mb X_h^K$. We show in
Lemma~6 that this matrix is asymptotically non-degenerate with
large probability when the design is regularly varying.

For technical reasons, we introduce a slightly different version
of the local polynomial estimator. We introduce a ``correction"
term in the matrix~$\mb X_h^K$.
\medskip

{\bf Definition 4.}
Given some $h > 0$, we consider $\wh f_h$ defined by (3.5)
with $\wh \tta_h$ the solution when it makes sense
(if $N_{n,h} = 0$ we take $\wh f_h = 0$) of the linear system
$$
\mbt X_h^K \tta = \mb Y_h^K,
\tag 3.10
$$
where
$$
\mbt X_h^K \eqdef \mb X_h^K + N_{n,h}^{1/2} \mb I_{k+1}
\ind{\lba(\mb X_h^K) \leq N_{n,h}^{1/2}},
$$
with $\lambda(M)$ being the smallest eigenvalue of a matrix~$M$
and $\mb I_{k+1}$ denoting the identity matrix in $\setR^{k+1}$.
\medskip

{\bf Remark.}
One can understand the definition of $\mbt X_h^K$ as follows: in
the ``good" case when $\mc X_h^K$ is non-degenerate in the
sense that its smallest eigenvalue is not too small, we solve
the system (3.7), while in the ``bad" case we still have a
control on the smallest eigenvalue of $\mbt X_h^K$,
since we always have $\lba(\mbt X_h^K) \geq N_{n,h}^{1/2}$.
\medskip

\qquad{\tensmc3.2. Bias-variance equilibrium}.
A main result on the local polynomial estimator is the
bias-variance decomposition. This is a classical result
presented many times in different forms: see
Cleveland~(1979), Goldenshluger and Nemirovski~(1997),
Korostelev and Tsybakov~(1993),
Spokoiny~(1998), Stone~(1980), Tsybakov~(1986, 2003).
The version in~Spokoiny~(1998) is close to the one presented here.
The differences are mostly related to the fact that the design
is random and that we consider a modified version of the local
polynomial estimator (see Definition~4).
We introduce the event
$$
\Omega_h^K \eqdef \{ X_1,\ldots,X_n \text{ are such that }
\lba(\mc X_h^K) > N_{n,h}^{-1/2} \text{ and } N_{n,h} > 0 \}.
\tag 3.11
$$
Note that on $\Omega_h^K$ the matrix $\mc X_h^K$ is invertible.
\medskip

{\bf Proposition 1} (Bias--variance decomposition).
{\it  Under Assumption~{\rm K} and if
$f \in \mc F_{h}(x_0,\omega)$, the following inequality holds on
the event $\Omega_h^K${\rm:}
$$
|\wh f_h(x_0) - f(x_0)| \leq
\lambda^{-1}(\mc X_h^K) \sqrt{k+1}
K_{\infty} \bigl( \omega(h) + \sigma N_{n,h}^{-1/2} |\gamma_h| \bigr),
\tag 3.12
$$
where $\gamma_h$ is, conditionally on $\mf X_n$, centered
Gaussian  such that $\Efm\{\gamma_h^2 \mid \mf X_n\}\leq1$.}
\medskip
\pagebreak

{\bf Remark.}
Inequality (3.12) holds conditionally on the design, on the
event $\Omega_h^K$. We will see that this event has a large
probability in the regular variation framework.
\medskip

\qquad{\tensmc3.3. Choice of the bandwidth}.
Now, like with any linear estimation procedure,
the problem is:
{\it how to choose the bandwidth $h$}?
In view of inequality (3.12) a natural bandwidth choice is
$$
H_n \eqdef \argmin_{h \in [0,1]}
\Big\{ \omega(h) \geq \frac{\sigma }{ \sqrt{N_{n,h}}} \Big\}.
\tag 3.13
$$
Such a bandwidth choice is well known, see, for
instance, Guerre~(2000). This choice stabilizes the procedure,
since it is sensitive to the design, which represents in the
model (1.1) the local amount of information. The estimator is
then defined by
$$
\wh f_n(x_0) \eqdef \wh f_{H_n}(x_0),
$$
where $\wh f_h$ is given by Definition~4 and $H_n$ is defined
by (3.13). The random bandwidth $H_n$ is close in probability
to the theoretical deterministic bandwidth $h_n$ defined by
(2.5) in view of the following proposition.
\medskip

{\bf Proposition 2.}
{\it Under Assumption~{\rm M} and if $\omega \in \RV(s)$ for any
$s > 0$, for any $0 < \von \leq 1/2$ there exists
$0 < \eta \leq \von$ such that
$$
\Pm\Big\{ \Big| \frac{H_n}{h_n} - 1 \Big| > \von \Big\}
\leq 4 \exp \Big(-\frac{\eta^2}{1+\eta/3} n F_{\nu}(h_n/2)\Bigl),
$$
where $F_{\nu}(h) \eqdef \int_0^h \nu(t)\, dt$.}
\medskip

If $n F_{\nu}(h_n / 2) \raro +\infty$ as $n \raro +\infty$
(this is the case when $\nu$ is regularly varying) this
inequality entails
$$
H_n = \big(1 + o_{\Pfm}(1)\big) h_n,
$$
where $o_{\Prob}(1)$ stands for a sequence going to~$0$
in probability under a probability~$\Prob$.

Proposition~3 motivates the regularly varying design choice.  It makes
a link between the behaviour of the counting process $N_{n,h}$ (that
appears in the variance term of (3.12)) and the behaviour of $\mu$
close to $x_0$.  Actually, the regular variation property (see
Definition~1) naturally appears under appropriate assumptions on the
asymptotic behaviour of $N_{n,h}$. Let us denote by $\Pm$ the joint
probability of the variables $(X_i)$.
\medskip

{\bf Proposition 3.}
{\it If Assumption~{\rm M} holds with $\nu$ monotone, then
the following properties are equivalent{\rm:}

{\rm(1)} $\nu$ is regularly varying of index
$\beta\geq-1${\rm;}

{\rm(2)} there exist sequences of positive numbers $(\lambda_n)$
and $(\gamma_n)$ such that $\lim_n \gamma_n = 0$,
$\liminf_n n \lambda_n^{-1} > 0$,
$\gamma_{n+1} \sim \gamma_{n}$ as $n \raro+\infty$
and a continuous function $\phi \: \setR^+ \raro \setR^+$
such that for any $C > 0${\rm:}
$$
\Em\{ N_{n,C \gamma_n} \} \sim \phi(C) \lambda_n
\qquad\text{as}\quad n\raro +\infty;
$$

{\rm(3)} there exist $(\lba_n)$, $(\gamma_n)$, and $\phi$ as before
such that for any $C > 0$ and $\von > 0${\rm:}
$$
\lim_{n \raro +\infty} \frac{n}{\lambda_n} \Pm\Bigl\{
\Bigl|\frac{N_{n,C\gamma_n}}{\phi(C)\lambda_n} - 1 \Bigr|
> \von\Bigr\} = 0.
$$}

The proof is delayed until Section~6. Mainly, it is a
consequence of the sequence characterization of regular
variation (see in the Appendix).
\medskip

{\bf 4. Upper Bounds for $\wh f_{H_n}(x_0)$}

\qquad{\tensmc4.1. Conditional on the design}.
When no assumptions on the behavior of the  design density are
made, we can work conditionally on the design. For $\lba > 0$ we
define the event
$$
\E_{\lba} \eqdef \{ \lba_n > \lba \},
$$
where $\lba_n \eqdef \lba(\mc X_{H_n}^K)$. Note that
$\E_{\lba}\in \mf X_n$. We also define the constant
$$
m(p)\eqdef\sqrt{2/\pi} \int_{\setR^+} (1 + t)^p
\exp(-t^2/2)\, dt.
$$

{\bf Proposition 4.}
{\it  Under Assumption~{\rm K}, if
$\lba$ is such that  
$\lba^2 N_{n,H_n} \geq 1$
and $n \geq k + 1$,
we have on $\E_{\lba}${\rm:}
$$
\sup_{f \in \mc F_{H_n}(x_0, \omega)}
\Efm\big\{ |\wh f_n(x_0) - f(x_0) |^p \mid \mf X_n \big\}
\leq m(p) \lba^{-p} K_{\infty}^p (k+1)^{p/2} R_n^p,
$$
where $R_n \eqdef \omega(H_n)$.}
\medskip

\qquad{\tensmc4.2. When the design is regularly varying}.
Proposition~5 below gives an upper bound for the estimator
$\wh f_{H_n}(x_0)$ when the design density is regularly varying.
This proposition can be viewed as a deterministic counterpart to
Proposition~4.

Let $\lba_{\beta, K}$ be the smallest eigenvalue of the symmetric
and positive matrix with entries, for $0 \leq j,l \leq k$:
$$
(\mc X_{\beta,K})_{j,l} = \frac{\beta+1}{2}
\bigl( 1+(-1)^{j+l} \bigr)\int_0^1 y^{j+l+\beta} K(y)\, dy.
\tag 4.1
$$
Note that in view of Lemma~6 we have $\lba_{\beta,K} >0$.
\medskip

{\bf Proposition 5.}
{\it Let $\varrho > 1$ and let $h_n$ be defined by {\rm(2.5)}.
Let $(\alpha_n)$ be a sequence of positive numbers such that
$\alpha_n = O(n^{\gamma})$ for some $\gamma > 0$.
If $\mu \in \mc R(x_0,\beta)$ with $\beta > -1$ and
$\omega \in \RV(s)$, we have for any $p > 0${\rm:}
$$
\limsup_n \sup_{{f \in \Sigma_{\varrho h_n,\alpha_n}(x_0, \omega)}}
\Efm\{ r_n^{-p} |\wh  f_n(x_0) -   f(x_0)|^p \}
\leq C \lba_{\beta,K}^{-p},
\tag 4.2
$$
where $r_n \eqdef \omega(h_n)$ satisfies
$$
r_n \sim \sigma^{2s/(1+2s+\beta)} n^{-s/(1+2s+\beta)}
\ell_{\omega,\nu}(1/n)
\qquad\text{as}\quad n \raro +\infty,
$$
with $\ell_{\omega,\nu}$ slowly varying and where
$C = 4^{s/(1+2s+\beta)}(k+1)^{p/2} m(p) K_{\infty}^p$.}
\medskip
\pagebreak

{\bf Remark.}
Under H{\"o}lder regularity with radius $r$ we have
$$
r_n \sim \sigma^{2s/(1+2s+\beta)} r^{(\beta+1)/(1+2s+\beta)}
n^{-s/(1+2s+\beta)} \ell_{s,\nu}(1/n)
\qquad\text{as}\quad n \raro+\infty.
$$
\medskip

{\bf 5. Discussion}

\qquad{\tensmc5.1. About Assumption~M}.
As stated previously, Assumption~M means that the
design distribution is symmetric around $x_0$ close to this
point.
When it is not the case, and if there are two functions
$\nu^{-}\in \RV(\beta^-)$, $\nu^{+} \in \RV(\beta^+)$
for $\beta^{-}, \beta^+ \geq-1$ and $\eta^{-}, \eta^{+} > 0$
such that for any $x \in [x_0-\eta^{-}, x_0 + \eta^{+}]$:
$$
\mu(x) = \nu^{+}(x - x_0) \ind{ x_0 \leq x \leq x_0
+
\eta^{+}} + \nu^{-}(x_0 - x) \ind{ x_0 - \eta^- \leq x < x_0},
$$
we can easily prove that the minimax convergence rate
is the fastest among the two possible ones, which is (2.4)
for the choice of $\beta = \beta^- \wedge \beta^+$.
To prove the upper bound we can use the same estimator as in
Section~3 with a non-symmetric choice of the bandwidth,
or more roughly we can ``throw away" the observations on the
side of $x_0$ corresponding to the largest index of regular
variation (when~$\mu$ is known).
\medskip

\qquad{\tensmc5.2. On Theorem~1 and Propositions~4 and~5}.
Since we are interested in the estimation of $f$ at $x_0$, we
need only a regularity assumption in some neighbourhood of this
point. Note that the minimax risks are computed over a class
where the regularity assumption holds in a decreasing interval
as~$n$ increases.

It appears that a natural choice of the size of this interval
is the theoretical bandwidth of estimation $h_n$, since it is
the minimum we need for the proof of the upper bounds. To state
an upper bound with the ``design-adaptive" estimator
$\wh f_{H_n}(x_0)$~--- in the sense
that it does not depend on the behavior of the design density
close to~$x_0$
(via the parameter $\beta$ for instance)~--- we need a
smoothness control in a slightly larger neighbourhood size
than~$h_n$ (see the parameter $\varrho$ in Proposition~5).

More precisely, to prove in Proposition~5 that~$r_n$ is an upper
bound, we use, in particular, Proposition~2 with
$\von = \varrho - 1$ in order to control the random bandwidth
$H_n$ by $h_n$. Thus, the parameter $\varrho$ is indispensable
for the proof of Proposition~5.
Note that we do not need such a parameter in Theorem~1 since we
use the estimator with the deterministic bandwidth $h_n$ to
prove the upper bound part of the theorem. Of course, this
estimator in unfeasible from a practical point of view since
$h_n$ heavily depends on $\mu$, which is hardly known in
practice. This is the reason why we state Proposition~5,
which tells us that the estimator with the data-driven bandwidth
$H_n$ converges with the same rate.
\medskip

\qquad{\tensmc5.3. On Theorem~2}.
In the $\Gamma$-variation framework, for the proof of the upper
bound part of Theorem~2 we use an estimator depending on~$\mu$.
Again, such an estimator is unfeasible from a practical point of
view. Anyway, this framework is considered only for theoretical
purposes, since from a practical point of view nothing can
be done in this case: there is no observations at the point of
estimation. This is precisely what Theorem~2 and the
corresponding example tell us, in the sense that the minimax
rate is very slow.
\medskip

\qquad{\tensmc5.4. About the $\Gamma$-varying design case}.
For the proof of the upper bound part in Theorem~2
we can consider an estimator different from the classical
regressogram (see the proof of the theorem). If $K$ is a
kernel satisfying Assumption~K, we define
$$
\wt f_n(x_0) \eqdef
\frac{ \sumin Y_i \big(K\big(\frac{X_i-h_n-x_0}{\rho(h_n)}\big)
+ K \big(\frac{X_i+h_n-x_0}{\rho(h_n)}\big)\big)}
{ \sumin K \big( \frac{ X_i - h_n - x_0 }{ \rho(h_n)} \big)
+ K \big( \frac{X_i + h_n - x_0}{ \rho(h_n) } \big) },
$$
where $h_n$ is defined by (2.5). The point is that since
$\supp K \subset [-1,1]$, this estimator makes a local average
of the observations $Y_i$ such that
$X_i \in [x_0 - h - \rho(h), x_0 - h + \rho(h)] \cup
[x_0 + h - \rho(h), x_0 + h + \rho(h)]$, which does not
contain the point of estimation $x_0$ for~$n$ large enough,
since $\lim_{h \raro 0^+} \rho(h) / h = 0$ (see Appendix).
In spite of this, we can prove that $\wt f_n(x_0)$ converges
with the rate $r_n$. We can understand this as follows: since
there is no information at $x_0$, the procedure actually
``catches" the information ``far" from~$x_0$.
This fact shows that again, the $\Gamma$-varying design is an
extreme case.
\medskip

\qquad{\tensmc5.5. More technical remarks}

$\bullet$ About Assumption~K, the first assumption is used to
make the kernel $K$ localize the information
around the point of estimation $x_0$ (see (3.2)).
The last one is technical and used in the proof of Lemma~6.
The two other ones are used for the sake of simplicity, since we
only really need the kernel to be bounded from above.

$\bullet$  When $\beta = -1$, Theorem~1 holds only for small
regularities $0 < s \leq 1$. For technical reasons, we
were not able to prove the upper bound when $s > 1$ and
$\beta =  -1$. More precisely, in this case we have $k=0$ and in
view of (3.4) it is clear that the local polynomial estimator is
a Nadaraya--Watson estimator defined by
$$
\wh f_n(x_0) = \frac{ \sumin Y_i K \bigl( \frac{X_i-x_0}{h_n}\bigr)}
{ \sumin K \bigl( \frac{X_i-x_0}{h_n} \bigr)}.
$$
When $s > 1$, we have to use a local polynomial estimator. The
problem is then in the asymptotic control of the smallest
eigenvalue  of $\mb X_{h_n}^K$ (see Lemma~6) and to do so we use
an average (Abelian) transform property of regularly varying
functions, which is (see Appendix):
$$
\lim_{h\ \raro 0^+} \frac{1}{\ell_{\nu}(h)}
\int y^{\alpha} K(y)\ell_{\nu}(yh)\, \frac{dy}{y}
= \cases
\int y^{\alpha-1} K(y)\,dy &\text{when } \alpha > 0, \\
+\infty                    &\text{when } \alpha = 0.
\endcases
$$
Thus the only way to have a limit for both cases is to assume
$K(y) = O(|y|^{\eta})$ for some $\eta > 0$, but the obtained
upper bound rate in this case would be slower than the lower
bound.
\medskip

{\bf 6. Proofs}

\qquad{\tensmc6.1. Proof of the main results}
\smallskip

{\it Proof of Theorem}~1.
First we prove the upper bound part of equation~(2.4)
when $\beta > -1$. We consider the estimator
$\wh f_n(x_0) = \wh f_{h_n}(x_0)$, where $\wh f_h$ is given
by Definition~4 with~$h_n$ given by equation (2.5),
and we define $r_n = \omega(h_n)$.
Let $0 < \von \leq \frac{1}{2}$.
We introduce the event
$$
\mc B_{n,\von} \eqdef \bigl\{ |\lba(\mc X_{h_n}^K) -
\lba_{\beta,K} | \leq \von \bigr\}
\cap \Bigl\{ \Bigl|\frac{N_{n,h_n}}{2 n F_{\nu}(h_n)} - 1 \Bigr|
\leq \von \Bigr\}.
$$
Since $\lim_n nF_{\nu}(h_n) = +\infty$ (see, for instance,
Lemma~4),
we have $\mc B_{n,\von} \subset \Omega_{h_n}^K$
for~$n$ large enough (see (3.11)) and, in particular,
on the event $\mc B_{n,\von}$
the matrix $\mb X_{h_n}^K$ is invertible.
Then using Proposition~1
and since $f \in \mc F_{h_n}(x_0,\omega)$,
we get:
$$
\align
| \wh f_n(x_0) - f(x_0) | \ind{\mc B_{n,\von}}
&\leq
(\lba_{\beta,K} - \von)^{-1} \sqrt{k+1} K_{\infty}
\Bigl(\omega(h_n)
+ \frac{\sigma}{\sqrt{(2 - \von) n F_{\nu}(h_n)}}|\gamma_{h_n}| \Bigr)
\\
&\leq
(\lba_{\beta,K} - \von)^{-1} \sqrt{k+1} K_{\infty}
\omega(h_n) (1 + |\gamma_{h_n}|),
\endalign
$$
where we last used the definition of $h_n$.
Since, conditionally on $\mf X_n$,
$\gamma_{h_n}$ is centered Gaussian such that
$\Efm\{  \gamma_{h_n}^2 \mid \mf X_n \} \leq 1$, we get for any $p > 0$:
$$
\sup_{f \in \mc F_{h_n}(x_0,\omega)}
\Efm\big\{ r_n^{-p} |\wh f_n(x_0) - f(x_0) |^p \,
\ind{\mc B_{n,\von}}\mid \mf X_n \big\}
\leq (\lba_{\beta,K} - \von)^{-p} (k+1)^{p/2} K_{\infty}^p m(p),
$$
where $m(p)$ is defined in Section~4.
Now we work on the complement $\mc B_{n,\von}^{c}$.
We use Lemmas~2 and~6 to control the probability of
$\mc B_{n,\von}$ and we recall that
$\alpha_n =  O(n^{\gamma})$ for some $\gamma > 0$.
When $N_{n,h_n} = 0$ we have  $\wh f_n(x_0) = 0$ by definition
and then
$$
\sup_{f \in \mc U(\alpha_n)}
\Efm\big\{ r_n^{-p} | \wh f_n(x_0) - f(x_0) |^p\,
\ind{\mc B_{n,\von}^{c}} \big\}
\leq (\alpha_n r_n^{-1})^p\,
\Pfm\{\mc B_{n,\von}^c \} = o_n(1).
$$
Then we assume $N_{n,h_n} > 0$. Using Lemma~3
we get:
$$
\align
&\sup_{f \in \mc U(\alpha_n)}
\Efm\big\{ r_n^{-p} |\wh f_n(x_0) - f(x_0)|^p\, \ind{\mc B_{n,\von}^c}
\big\}
\\
&\qquad
\leq 2^p r_n^{-p}
\Big(\sqrt{\Efm\{ |\wh f_n(x_0)|^{2p} \}} + \alpha_n^{p}\Big)
\sqrt{\Pm\{ \mc B_{n,\von}^c \}}
\\
&\qquad
\leq 2^p (\alpha_n r_n^{-1})^p
\big( \sqrt{n^p C_{\sigma, k, 2p}} + 1\big)
\sqrt{\Pm\{ \mc B_{n,\von}^c \}} = o_n(1),
\endalign
$$
and thus we have proved that $r_n$ is an upper bound of the
minimax risk (2.4) when $\beta > -1$.

When $\beta = -1$ and $0 < s \leq 1$, we have $k = 0$ and the
matrix $\mc X_{h_n}^K$ is $1 \times 1$ sized and equal to
$\overline{K}_{n,h_n,0}$ (see equation (6.5)).
The bias--variance equation (3.12) becomes in this case:
$$
|\wh f_n(x_0) - f(x_0)| \leq (\overline{K}_{n,h_n,0})^{-1}
K_{\infty}\big(\omega(h_n)
+ \sigma N_{n,h_n}^{-1/2} |\gamma_{h_n}| \big).
$$
Consider the event
$$
\mc C_{n,\von} = \Big\{ \Big| \frac{N_{n,h_n}}{2n F_{\nu}(h_n)}
-1 \Big| \leq \von \Big\}
\cap
\Bigl\{ \Bigl| \frac{K_{n,h_n,0}}{2n F_{\nu}(h_n)} - K(0) \Bigr|
\leq \von \Bigr\}.
$$
We note that the probability of $\mc C_{n,\von}$ is controlled
by Lemma~2 and equation (6.8) in Lemma~5.
Then we can proceed as previously to prove that $r_n$ is an
upper bound when $\beta = -1$ and we have proved that $r_n$ is
an upper bound for the left-hand side of (2.4).
Using Proposition~6 we also have that $r_n$ is a lower bound for
the left part of (2.4). The conclusion follows from Lemma~4.
\qed
\medskip

{\it Proof of Theorem}~2.
The proof is similar to that of Theorem~1.
For the proof of the upper bound part in (2.7)
we use the regressogram estimator defined by
$$
\wh f_n(x_0) \eqdef \cases
\dsize  \frac{\sumin Y_i \ind{|X_i - x_0| \leq h_n}}
{N_{n,h_n}} &\text{if}\quad N_{n,h_n} > 0,
\\
0           &\text{if}\quad N_{n,h_n} = 0.
\endcases
$$
Let $0 < \von \leq 1/2$. On the event
$
\mc D_{n,\von} \eqdef \Big \{ \Big|
\frac{N_{n, h_n}}{2 n F_{\nu}(h_n)} - 1 \Big| \leq \von \Big \}$
we clearly have $N_{n,h_n} > 0$ and since
$f \in \mc F_{h_n}(x_0,\omega)$, we have
$$
|\wh f_n(x_0) - f(x_0)|
\leq \omega(h_n) + \sigma N_{n,h_n}^{-1/2}|v_{n}|
\leq \omega(h_n) (1 - \von)^{-1/2} (1 + |v_n|),
$$
where $v_n \eqdef \frac{1}{\sigma \sqrt{N_{n,h_n}}} \sumin \xi_i
\ind{|X_i - x_0| \leq h_n}$ is, conditionally on $\mf X_n$,
standard Gaussian. Then we get
$$
\sup_{f \in \mc F_{h_n}(x_0, \omega)}
\Efm\big\{ |\wh f_n(x_0) - f(x_0)|^p \ind{\mc D_{n,\von}} \big\}
\leq r_n^p(1 - \von)^{-p/2} m(p).
$$
Now we work on $\mc D_{n,\von}^c$. If $N_{n,h_n} = 0$,
we get using Lemma~2 and since $\alpha_n = O(r_n^{-\gamma})$:
$$
\align
\sup_{f \in \mc U(\alpha_n) }
\Efm\big\{ | \wh f_n(x_0) - f(x_0)|^p\,
\ind{\mc D_{n,\von}^c} \big\}
&\leq \alpha_n^p \Pm\{ \mc D_{n,\von}^c \}
\\
&= O( r_n^{-\gamma p}) \exp \Big(
-\frac{ \von^2 \sigma^2 }{ 1 + \von / 3} r_n^{-2} \Big)
= o_n(1),
\endalign
$$
since $\alpha_n = O(r_n^{-\gamma})$. If $N_{n,h_n} > 0$,
since $|\wh f_n(x_0)| \leq \alpha_n + \sigma |v_n|$, we get
$$
\sup_{f \in \mc U(\alpha_n) } \Efm\big\{ | \wh f_n(x_0) - f(x_0)|^p
\ind{\mc D_{n,\von}^c} \big\}
\leq 2^p \alpha_n^p (1 +\sqrt{C_{\sigma, 0, p}} )
\sqrt{ \Pm\{ \mc D_{n,\von}^c \} } = o_n(1),
$$
where $C_{\sigma, 0, p}$ is the same as in the proof of
Theorem~1. Thus we have proved that $r_n$ is an upper bound.
The lower bound is given by Proposition~6,
and the conclusion follows from Lemma~4.
\qed
\medskip

In the sequel, $\prodsca{\cdot}{\cdot}$ denotes the Euclidean
scalar product on $\setR^{k+1}$, \
$e_1 =\mathbreak  (1,0,\ldots,0) \in \setR^{k+1}$, \
$\norminfty{\cdot}$ stands for the sup norm in $\setR^{k+1}$,
and $\norm{\cdot}$ stands for the Euclidean norm in
$\setR^{k+1}$.
\medskip

{\it Proof of Proposition}~1.
On $\Omega_h^K$ we have in view of Definition~4
that $\mbt X_h^K = \mb X_h^K$ and $\mb X_h^K$ is invertible. Let
$0 < \von \leq 1/2$ and $n \geq 1$. We can find a polynomial
$P_f^{n,\von}$ of order $k$ such that
$$
\sup_{|x-x_0| \leq h} |f(x) - P_f^{n,\von}(x)|
\leq
\inf_{P \in  \mc P_k} \sup_{|x - x_0| \leq h } |f(x) - P(x - x_0)|
+  \frac{\von}{\sqrt{n}}.
$$
In particular, with $h = 0$ we get $|f(x_0) - P_f^{n,\von}(x_0)|
\leq \frac{\von}{\sqrt{n}}$. Defining $\tta_h \in \setR^{k+1}$
such that $P_f^{n,\von} = P_{\tta_h}$ (see (3.6)) we get
$$
|\wh f_h(x_0) - f(x_0)| \leq \frac{\von}{\sqrt{n}}
+ |\prodsca{\wh \tta_h - \tta_h}{e_1}|
= \frac{\von}{\sqrt{n}}
+ |\prodsca{(\mb X_h^{K})^{-1}\mb X_h^K (\wh\tta_h - \tta_h)}{e_1}|.
$$
Then we have for $j \in \{0,\ldots,k \}$ by (3.4) and (1.1):
$$
\align
(\mb X_h^K (\wh \tta_h - \tta_h))_j
&= \prodscahk{\wh f_{h} - P_f^{n,\von}}{\phi_{j,h}}
= \prodscahk{Y - P_f^{n,\von}}{\phi_{j,h}} \\
&= \prodscahk{f - P_f^{n,\von}}{\phi_{j,h}}
+ \prodscahk{Y - f}{\phi_{j,h}} \\
&= \prodscahk{f - P_f^{n,\von}}{\phi_{j,h}}
+ \prodscahk{\xi}{\phi_{j,h}}
\eqdef B_{h,j} + V_{h,j},
\endalign
$$
thus $\mb X_h^K (\wh \tta_h - \tta_h) = B_h + V_h$. In view of
Assumption~K and since $f \in \mc  F_h(x_0,\omega)$, we have:
$$
|B_{h,j}| = |\prodscahk{f - P_f^{n,\von}}{\phi_{j,h}}|
\leq
\normhk{f - P_f^{n,\von}}\, \normhk{\phi_{j,h}}
\leq N_{n,h} K_{\infty} \Big(\omega(h) + \frac{\von}{\sqrt{n}}\Big),
$$
thus $\norminfty{B_h} \leq N_{n,h} K_{\infty}
(\omega(h) +\frac{\von}{\sqrt{n}})$.
Moreover, since $\lba^{-1}(\mc X_h) \leq  N_{n,h}^{1/2}
\leq n^{1/2}$ on $\Omega_{h,K}$, we have:
$$
\align
|\prodsca{(\mb X_h^K)^{-1}B_h}{e_1}|
&\leq \norm{(\mb X_h^K)^{-1}}\,\norm{B_h}
\leq \norm{(\mb X_h^K)^{-1}} \sqrt{k+1} \norminfty{B_h}
\\
&\leq \lba^{-1}(\mc X_h^K) \sqrt{k+1} K_{\infty} \omega(h)
+  \sqrt{k+1} K_{\infty} \von,
\endalign
$$
where we last used the fact that $\norm{M^{-1}} = \lba^{-1}(M)$
for a positive symmetric matrix. The variance term $V_h$ is
clearly, conditionally on $\Xn$, a centered Gaussian vector, and
its covariance matrix is equal to $\sigma^2 \mb X_h^{K^2}$.
Thus the random variable $\prodscahk{(\mb X_h^K)^{-1} V_h}{e_1}$
is, conditionally on $\Xn$, centered Gaussian of variance:
$$
\align
v_h^2 &= \sigma^2 \prodsca{e_1}{(\mb X_h^K)^{-1} \mb X_h^{K^2}
(\mb X_h^K)^{-1} e_1}
\leq \sigma^2 \prodsca{e_1}{(\mb X_h^K)^{-1}
\mb X_h^{K} (\mb X_h^K)^{-1} e_1}
\\
&= \sigma^2 \prodsca{e_1}{(\mb X_h^K)^{-1} e_1}
\leq \sigma^2 \norm{(\mb X_h^K)^{-1}}
= \sigma^2 N_{n,h}^{-1} \lba^{-1}(\mc X_h^K),
\endalign
$$
since $K \leq 1$. Then
$
\lba(\mc X_h^K) = \inf_{\norm{x} = 1} \prodsca{x}{\mc X_h^K x}
\leq \norm{\mc X_h^K e_1} \leq \sqrt{k+1}$,
since $\mc X_h^K$ is symmetric and its entries are smaller
than~$1$ in absolute value. Thus
$$
v_h^2 \leq \sigma^2 N_{n,h}^{-1} \lba^{-1}(\mc X_h^K)
\leq\sigma^2 N_{n,h}^{-1} (k+1) \lba^{-2}(\mc X_h^K),
$$
and the proposition follows.
\qed
\medskip

{\it Proof of Proposition}~2.
The proposition is a direct consequence of Lemmas~1 and~2.
\qed
\medskip

{\it Proof of Proposition}~3.
$(2) \Rightarrow (1)$: In view of Assumption~M one has for $n$
large enough
$$
\Em\{ N_{n,C\gamma_n} \} = 2 n \int_0^{C \gamma_n} \nu(x)\, dx
= 2 n F_{\nu}(C \gamma_n),
$$
thus $(2)$ entails
$2 n \lambda_n^{-1} F_{\nu}(C\gamma_n) \sim \phi(C)$ as
$n \raro +\infty$ and then $F_{\nu} \in \RV(\alpha)$ in view of
the   characterization (A.8) of regular variation.
Since $F_{\nu}(0) = 0$, we have more precisely
$F_{\nu}  \in \RV(\alpha)$ for $\alpha \geq 0$ and since $\nu$
is monotone, we  have $\nu \in \RV(\alpha-1)$ (see Appendix).

$(3) \Rightarrow (2)$: Let $\von > 0$. We define the event
$$
A_n(C, \von)
= \Bigl\{\Bigl| \frac{N_{n,C\gamma_n}}{\phi(C)\lambda_n} -1 \Bigr|
\leq \von \Bigr\}.
$$
Then:
$$
\align
\lambda_n^{-1} \Em\{ N_{n,C\gamma_n} \}
&= \lambda_n^{-1} \Em \bigl\{ N_{n,C\gamma_n} ( \ind{A_n(C,\von)}
+ \ind{A_n^{c}(C,\von)})\bigr\}\\
&\leq (1+\von) \phi(C)
+ n \lambda_n^{-1} \Pm\bigl\{ A_n^{c}(C, \von)\bigr\},
\endalign
$$
and then $\limsup_n \lambda_n^{-1} \Em\{ N_{n,C\gamma_n} \}
\leq  (1+\von) \phi(C)$. On the other hand,
$$
\lambda_n^{-1} \Em \{ N_{n,C\gamma_n} \}
\geq \lambda_n^{-1} \Em \{ N_{n,C\gamma_n} \ind{A_n(C,\von)} \}
\geq (1-\von) \phi(C) \Pm\{ A_n(C,\von) \},
$$
and then $\liminf_n \lambda_n^{-1} \Em\{ N_{n,C\gamma_n}\}
\geq (1-\von) \phi(C)$.

$(1) \Rightarrow (3)$: Let $\nu \in \RV(\beta)$ and
$0 < \von \leq  1/2$. If $\beta > -1$, we have
$F_{\nu} \in \RV(\beta+1)$ (see in the Appendix),
thus we can write $F_{\nu}(h) = h^{\beta+1} \ell_{F}(h)$,
where $\ell_{F}$ is slowly varying. We define
$\gamma_n =  n^{-1/(2(\beta+1))}$ when $\beta > -1$
and $\gamma_n = n^{-1}$ if  $\beta = -1$.
When $\beta = -1$, we have $F_{\nu} \in \RV(0)$
(see Appendix). We note that in both cases we have
$\lim_n \gamma_n = 0$  and $\gamma_{n+1} \sim \gamma_{n}$
as $n \raro +\infty$. In view of Lemma~2 we get for $n$
large enough
$$
\Pm \Bigl\{\Bigl|\frac{N_{n,C\gamma_n}}{\phi(C)\lba_n}-1\Bigr|
> \von \Bigr\}
\leq 2 \exp \Bigl( -\frac{\von^2}{1 + \von/3} \phi(C) \lba_n\Bigr),
$$
where we used the fact that $\ell_{F}$ is slowly varying and
where  we defined $\lba_n \eqdef 2 n F_{\nu}(\gamma_n)$ and
$\phi(C) \eqdef  C^{\beta+1}$. Then we clearly have
$\lim_n n \lba_n^{-1} = +\infty$
and the proposition follows.
\qed
\medskip

\qquad{\tensmc6.2. Proof of the upper bounds for $\wh f_{H_n}(x_0)$}
\smallskip

{\it Proof of Proposition}~4.
Since $\E_{\lba} \subset \Omega_{H_n}^K$,
(3.13) and Proposition~1 entail that uniformly in
$f \in \mc F_{H_n}(x_0, \omega)$ we have
$$
|\wh f_n(x_0) - f(x_0)| \leq \lba^{-1} \sqrt{k+1}
 K_{\infty} R_n (1 + |\gamma_{H_n}|),
$$
where $\gamma_{H_n}$ is, conditionally on $\mf X_n$, centered
Gaussian  such that $\Efm\{ \gamma_{H_n}^2 \mid \mf X_n \} \leq 1$.
The result follows by integration with respect to
$\Pfm(\cdot\mid\mf X_n)$.
\qed
\medskip

{\it Proof of Proposition}~5.
Let us define $\von \eqdef \varrho - 1$. We can assume without
loss of generality that $\von < \frac{1}{2} \wedge \lba_{\beta,K}$.
We consider the event $\mc A_{n,\von}$ from Lemma~6.
In view of this lemma we have $\mc A_{n,\von} \subset
\E_{\lba_{\beta,K} - \von} \cap \{ (1 - \von) h_n \leq
H_n \leq (1 + \von) h_n \}$ and then
$\mc F_{\varrho h_n}(x_0, \omega) \subset \mc F_{H_n}(x_0,\omega)$.
Thus using Proposition~4 we get
$$
\align
&\sup_{f \in \mc F_{\varrho h_n}(x_0,\omega)}
\Efm\{ |\wh f_n(x_0) - f(x_0)|^p \ind{\mc A_{n,\von}} \mid \mf X_n \}
\\
&\qquad
\leq m(p) (\lba_{\beta,K} - \von)^{-p} K_{\infty}^p (k+1)^{p/2} R_n^p
\\
&\qquad
\leq m(p) (\lba_{\beta,K} - \von)^{-p} K_{\infty}^p (k+1)^{p/2}
(1+\von)^{p(s+1)} r_n^p,
\endalign
$$
where we used equation (6.1) in the same way as in the proof of
Lemma~1 to obtain on $\mc A_{n,\von}$ that
$\omega(H_n) \leq (1 + \von)^{s+1} \omega(h_n)$.
On the complementary event $\mc A_{n,\von}^c$,
using inequality (6.11) and Lemma~3
and since $\alpha_n = O(n^{\gamma})$ for some $\gamma > 0$,
we get
$$
\align
&\sup_{f \in \mc U(\alpha_n)} \Efm\{ r_n^{-p}
|\wh f_n(x_0) - f(x_0)|^p \ind{\mc A_{n,\von}^c} \}
\\
&\qquad
\leq 2^p (\alpha_n r_n^{-1})^p ( \sqrt{n^p C_{\sigma, k, 2p} } + 1 )
\sqrt{\Pm\{ \mc A_{n,\von}^c \}} = o_n(1),
\endalign
$$
and (4.2) follows. The equivalent of $r_n$ is given by Lemma~4.
\qed
\medskip

\qquad{\tensmc6.3. Lemmas for the proof of the upper bounds}
\smallskip

{\bf Lemma 1.}
{\it If $\omega \in \RV(s)$ for any $s > 0$, then for any
$0 < \von \leq \frac{1}{2}$ there exists $0 < \eta \leq \von$
such that
$$
\Bigl\{ \Bigl| \frac{N_{n,(1-\von)h_n}}{2n F_{\nu}((1 - \von)h_n)}
- 1 \Bigr| \leq \eta \Bigr\}
\cap
\Bigl\{ \Bigl| \frac{N_{n, (1 + \von) h_n}}{2n F_{\nu}((1 + \von)h_n)}
- 1 \Bigr| \leq \eta \Bigr\}
\subset \Big\{ \Big| \frac{H_n}{h_n} - 1 \Big| \leq \von \Big\}.
$$}

{\it Proof}.
In view of (3.13) we have
$\{ H_n \leq (1 + \von) h_n \} = \{ N_{n,(1 + \von) h_n}
\geq  \sigma^2 \omega^{-2}((1 + \von)h_n) \}$.
Define $\von_1 \eqdef 1 - (1 - \von^2)^{-2} (1 + \von)^{-2s}$.
For $\von$ small enough, it is clear that $\von_1 > 0$.
We recall that $\ell_{\omega}$ stands for the slowly varying
term of $\omega$ (see Definition~2). Since (A.1) holds uniformly
on each compact set in $(0,+\infty)$,
we have for $n$ large enough that for any
$y\in[\frac{1}{2},\frac{3}{2}]$:
$$
(1 - \von^2) \ell_{\omega}(h_n) \leq \ell_{\omega}(y h_n)
\leq (1 + \von^2) \ell_{\omega}(h_n),
\tag 6.1
$$
so using (6.1) with $y = 1 + \von$ \
($\von \leq \frac{1}{2}$), we obtain in view of (2.5):
$$
\align
2 (1 - \von_1) n F_{\nu}((1 + \von)h_n)
&\geq (1 - \von^2)^{-2} (1 + \von)^{-2s} \sigma^2 \omega^{-2}(h_n)
\\
&= \sigma^2 \big( (1 + \von)h_n \big)^{-2s}
(1 - \von^2)^{-2} \ell_{\omega}^{-2}(h_n)
\\
&\geq \sigma^2 \omega((1 + \von) h_n)^{-2},
\endalign
$$
and then
$$
\{N_{n,(1 + \von) h_n}\geq 2(1-\von_1)nF_{\nu}((1 + \von) h_n)\}
\subset \{ H_n \leq (1 + \von) h_n \}.
$$
Using again (6.1) with $y = 1 - \von$ we get in the same way
$$
\{ N_{n,(1 - \von) h_n} < 2(1 + \von_1) n F_{\nu}((1 - \von)h_n)\}
\subset \{ H_n > (1 - \von) h_n \},
$$
and then
$$
\Bigl\{ \Bigl| \frac{N_{n, (1 - \von) h_n}}
{2n F_{\nu}((1 - \von)h_n)} - 1 \Bigr| \leq \von_1 \Bigr\}
\cap
\Bigl\{ \Bigl| \frac{N_{n, (1 + \von) h_n}}
{2n F_{\nu}((1 + \von)h_n)} - 1 \Bigr| \leq \von_1 \Bigr\}
\subset
\Big\{ \Big| \frac{H_n}{h_n} - 1 \Big| \leq \von \Big\}.
$$
Now the result follows for the choice $\eta = \von \wedge \von_1$.
\qed
\medskip

{\bf Lemma 2.}
{\it Under Assumption~{\rm M},
we have for any $\von,h>0${\rm:}
$$
\Pm \Bigl\{\Bigl| \frac{N_{n, h}} {2n F_{\nu}(h)}
- 1 \Bigr| > \von \Bigr\}
\leq 2 \exp \Bigl( - \frac{\von^2}{1 + \von / 3} n F_{\nu}(h) \Bigr).
$$}

{\it Proof}.
It suffices to apply the Bernstein inequality to the sum of
independent random variables
$Z_i = \ind{|X_i - x_0| \leq h} - \Pm\{|X_1 - x_0| \leq h\}$
for $i=1,\ldots,n$.
\qed
\medskip

{\bf Lemma 3.}
{\it For any $p>0$ and $h>0$ the estimator $\wh f_h$
{\rm(}see Definition~{\rm4)} satisfies
$$
\sup_{f \in \mc U(\alpha)}
\Efm\{|\wh f_h(x_0)|^p \mid\mf X_n \}
\leq  C_{\sigma, k, p} (\alpha \sqrt{n})^p,
$$
where $C_{\sigma, k, p} \eqdef (k+1)^{p/2} \sqrt{ 2 / \pi}
\int_{\setR^+} (1 + \sigma t)^p \exp(-t^2 / 2)\, dt$.}
\medskip

{\it Proof}.
When $N_{n,h} = 0$, we have $\wh f_h = 0$ by definition
and the result is obvious, so we assume $N_{n,h} > 0$.
Using the fact that $\lba(A+B) \geq \lba(A) + \lba(B)$
when $A$ and $B$ are symmetric non-negative matrices
we get $\lba(\mbt X_h^K) \geq  N_{n,h}^{1/2} > 0$,
thus $\mbt X_h^K$ is invertible. Equation~(3.10)
entails $|\wh f_h(x_0)|
= |\prodsca{(\mbt X_h^K)^{-1} \mbt X_h^K \wh \tta_h}{e_1} |
= |\prodsca{(\mbt X_h^K)^{-1} \mb Y_h }{e_1} |$.
In view of (1.1) we can decompose for $j \in \{ 0,\ldots,k \}$:
$$
(\mb Y_h)_j = \prodscahk{Y}{\phi_{j,h}}
= \prodscahk{f}{\phi_{j,h}} + \prodscahk{\xi}{\phi_{j,h}}
\eqdef  B_{h,j} + V_{h,j}.
$$
Since $f \in \mc U(\alpha)$, we have under Assumption~K
that $|B_{h,j}| \leq \alpha N_{n,h}$,
thus $\norminfty{B_h} \leq \alpha N_{n,h}$.
As in the proof of Proposition~1 we have that
$\prodsca{(\mbt{X}_h^K)^{-1} V_h}{e_1}$ is,
conditionally on $\mf X_n$,
centered Gaussian with variance
$$
\align
v_h^2 &= \sigma^2 \prodsca{e_1}{(\mbt X_h^K)^{-1} \mb X_h^{K^2}
(\mbt X_h^K)^{-1} e_1} \\
&\leq \sigma^2 \prodsca{e_1}{(\mbt
X_h^K)^{-1} \mb X_h^{K} (\mbt X_h^K)^{-1} e_1}
\leq \sigma^2 \norm{(\mbt X_h^K)^{-1}}^2 \norm{\mb X_h^{K}}.
\endalign
$$
Assumption~K entails that all the elements
of the matrix $\mb X_h^K$ are smaller than $N_{n,h}$, thus
$\norm{\mb X_h^K} \leq (k+1) N_{n,h}$. Since $\mbt X_h^K$ is
symmetric, we get
$\norm{(\mbt X_h^K)^{-1}} = \lba^{-1}(\mbt X_h^K)
\leq N_{n,h}^{-1/2}$, and then $v_h^2 \leq \sigma^2 (k+1)$.
Finally, we have
$$
\align
|\wh f_h(x_0)|
&\leq |\prodsca{(\mbt{X}_h^K)^{-1} B_h}{e_1}|
+    |\prodsca{ (\mbt{X}_h^K)^{-1} V_h}{e_1}|
\\
&\leq
\norm{(\mbt{X}_h^K)^{-1}}\, \norm{B_h}
+ \sigma \sqrt{k+1} |\gamma_h|
\leq
\sqrt{k+1} (\alpha \sqrt{n} + \sigma |\gamma_h|),
\endalign
$$
where $\gamma_h$ is, conditionally on $\mf X_n$,
centered Gaussian with variance smaller than~$1$.
The result follows by integrating with respect to
$\Pfm( \cdot \mid \mf X_n )$.
\qed
\medskip

{\bf Lemma 4.}
{\it If $\nu \in \RV(\beta)$, $\omega \in \RV(s)$ for $s>0$
and the sequence $(h_n)$ is defined by {\rm(2.5)}
then the rate $r_n  = \omega(h_n)$ satisfies
$$
r_n \sim c_{s,\beta} \sigma^{2s/(1+2s+\beta)} n^{-s/(1+2s+\beta)}
\ell_{\omega,\nu}(1/n)\qquad\text{as}\quad n\raro +\infty,
\tag 6.2
$$
where $\ell_{\omega,\nu}$ is slowly varying and $c_{s,\beta} =
4^{s/(1+2s+\beta)}$. When $\omega(h) = r h^s$
{\rm(}H{\"o}lder regularity{\rm)} for $r>0$,
we have more precisely{\rm:}
$$
r_n \sim c_{s,\beta} \sigma^{2s/(1+2s+\beta)}
r^{(\beta+1) /(1+2s+\beta)} n^{-s/(1+2s+\beta)}
\ell_{s,\nu}(1/n)
\quad\text{as}\quad n \raro +\infty,
\tag 6.3
$$
where $\ell_{s,\nu}$ is slowly varying.
It is noteworthy that when $\beta=-1$ the result becomes{\rm:}
$$
r_n \sim 2 \sigma n^{-1/2} \ell_{\omega,\nu}(1/n)
\qquad\text{as}\quad n \raro +\infty.
$$
When $\nu \in \GV(\rho)$, we have
$$
r_n \sim \ell_{\omega,\nu}(1 / n),
\tag 6.4
$$
where $\ell_{\omega,\nu}$ is slowly varying.}
\medskip

{\it Proof}.
Denote $F_{\nu}(h) \eqdef \int_0^h \nu(t)\, dt$
and let $G(h) =  \omega^2(h) F_{\nu}(h)$.
When $\beta > -1$, we have $F_{\nu} \in \RV(\beta+1)$
(see the Appendix) and when $\beta=-1$,
$F_{\nu}$ is slowly varying.
Thus $G \in \RV(1+2s+\beta)$ for any $\beta \geq -1$.
The function $G$ is continuous and such that
$\lim_{h \raro 0^+} G(h) = 0$ in view of (A.2), since
$1+2s+\beta > 0$. Then, for $n$ large enough,   
$h_n = G^{\laro}(\sigma^2/ (4n))$, where
$G^{\laro}(h) \eqdef \inf\{ y \geq 0 \mid G(y) \geq h \}$
is the generalized inverse of $G$. Then in view of
(A.8) we have $G^{\laro} \in \RV(1/(1+2s+\beta))$
and then $\omega \circ G^{\laro} \in \RV(s/(1+2s+\beta))$
(see Appendix). Thus we can write
$\omega \circ  G^{\laro}(h) = h^{s/(1+2s+\beta)}
\ell_{\omega,\nu}(h)$, where $\ell_{\omega,\nu}$ is a slowly
varying function. Thus:
$$
\align
r_n &= \omega \Big( G^{\laro} \Big( \frac{\sigma^2}{4n} \Big) \Big)
= c_{s,\beta} \sigma^{2s/(1+2s+\beta)} n^{-s/(1+2s+\beta)}
\ell_{\omega,\nu} \Big (\frac{\sigma^2}{4 n} \Big)
\\
&\sim c_{s,\beta} \sigma^{2s/(1+2s+\beta)} n^{-s/(1+2s+\beta)}
\ell_{\omega,\nu}(1/n)
\qquad\text{as}\quad n \raro +\infty,
\endalign
$$
since $\ell$ is slowly varying. When $\omega(h) = r h^s$,
we can write more precisely
$h_n = G^{\laro}(\sigma^2 / (4 r^2 n))$,
where $G(h) = h^{2s} F_{\nu}(h)$, so (6.2) and (6.3) follow.

Let $y \in \setR$. Using (A.9) and the uniformity in (A.1)
we get
$\lim_{h \raro 0^+} \ellom(h + y \rho(h)) / \ellom(h) = 1$,
thus $\lim_{h \raro 0^+} \omega(h + y \rho(h))/\omega(h) = 1$.
Moreover,
since $\GV(\rho)$ is stable under integration
(see Appendix)
we have $F_{\nu} \in \GV(\rho)$, thus
$\lim_{h \raro 0^+} G(h + y \rho(y))/G(h) = \exp(y)$ and then
$G \in \GV(\rho)$. For $n$ large  enough, $h_n$ is well defined
and given by $h_n = G^{\laro}(\sigma^2/ (4n))$.
Since $G^{\laro} \in \PV(\ell)$ for
$\ell = \rho \circ \nu^{\laro} \in \RV(0)$
(see Appendix),
$G^{\laro}$ belongs, in particular, to $\RV(0)$ in view
of~(A.11) and then $r_n = \omega \circ G^{\laro}(\sigma^2/(4 n))$,
where $\omega \circ G^{\laro} \in \RV(0)$.
Thus $r_n \sim \omega \circ G^{\laro}(n^{-1})$
as $n \raro +\infty$ and (6.4) follows with
$\ell_{\omega,\nu} = \omega \circ G^{\laro}$.
\qed
\medskip

\qquad{\tensmc Study of the terms $\lambda(\mc X_{h_n}^{K})$
and $\lambda(\mc X_{H_n}^{K})$}.
We recall that the matrix $\mc X_{h,K}$ is defined as the
symmetric and non-negative matrix with entries
$(\mc X_{h,K})_{j,l} = \overline{K}_{n,h,j+l}$
for $0 \leq j,l \leq k$, where:
$$
\overline{K}_{n,h,\alpha} \eqdef \frac{1}{N_{n,h}} \sumin
\Bigl(\frac{X_i-x_0}{h}\Bigr)^{\alpha}
K\Bigl(\frac{X_i-x_0}{h}\Bigr),
\tag 6.5
$$
for $\alpha \in \setN$. Define
$K_{n,h,\alpha} \eqdef N_{n,h}\overline{K}_{n,h,\alpha}$ and
$$
K_{\alpha,\beta} \eqdef (1 + (-1)^{\alpha})
\int_0^1 y^{\alpha+\beta} K(y)\, dy.
\tag 6.6
$$
We define for any $\von > 0$ the event
$$
\D_{n,h,\alpha,K,\von} \eqdef \Bigl\{ \Bigl|
\frac{K_{n, h, \alpha}}{n F_{\nu}(h)}
- (\beta + 1) K_{\alpha,\beta}\Bigr| \leq \von \Bigr\}.
$$

{\bf Lemma 5.}
{\it Let $\alpha \in \setN$ and $\von > 0$. Under
Assumption~{\rm K} and if $\mu \in \mc R(x_0,\beta)$ with
$\beta>-1$, then for any positive sequence $(\gamma_n)$
going to~$0$ we have for $n$ large enough
$$
\Pm \bigl\{ \D_{n,\gamma_n,\alpha,K,\von}^c \bigr\} \leq 2
\exp \Bigl(-\frac{\von^2}{8(2 + \von/3)} n F_{\nu}(\gamma_n)\Bigr).
\tag 6.7
$$
When $\beta = -1$ we have{\rm:}
$$
\Pm \Bigl\{ \Bigl| \frac{K_{n,\gamma_n,0}}{n F_{\nu}(\gamma_n)}
-  2 K(0) \Bigr| > \von \Bigr\}
\leq
2 \exp \Bigl( - \frac{ \von^2 }{ 8 (2 + \von/3) }
n F_{\nu}(\gamma_n) \Bigr).
\tag 6.8
$$}

{\it Proof}.
First we prove (6.7). We define
$Q_{i,n,\alpha}\eqdef
\bigl( \frac{X_i-x_0}{\gamma_n} \bigr)^{\alpha}
K\bigl(\frac{X_i-x_0}{\gamma_n}\bigr)$,
$Z_{i,n,\alpha} \eqdef
Q_{i,n,\alpha} - \Em\{ Q_{i,n,\alpha} \}$.
Since $\mu \in \mc  R(x_0,\beta)$, one has for $i=1,\ldots,n$:
$$
\frac{1}{n F_{\nu}(\gamma_n)} \Em \{ Q_{i,n,\alpha} \}
= \frac{\gamma_n \nu(\gamma_n)}{F_{\nu}(\gamma_n) }
\frac{1 + (-1)^{\alpha}}{\ell_{\nu}(\gamma_n)}
\int_0^1 y^{\alpha+\beta} K(y) \ell_{\nu}( y \gamma_n)\, dy,
$$
where we used Assumption~K and the fact that
$[x_0-\gamma_n,x_0+\gamma_n] \subset W$
for $n$ large enough.
Then equations (A.3) and (A.4) entail:
$$
\lim_n \frac{1}{n F_{\nu}(\gamma_n)} \Em \{ Q_{i,n,\alpha} \}
= (\beta+1) K_{\alpha,\beta},
$$
and for $n$ large enough:
$$
\D_{n,\gamma_n,\alpha,K,\von}^c \subset
\biggl\{ \Bigl|\frac{1}{ n F_{\nu}(\gamma_n)}
\sum_{i=1}^n Z_{i,n,\alpha}\Bigr| > \von / 2 \biggr\}.
\tag 6.9
$$
In view of Assumption~K we have
$\Em\{{Z_{i,n,\alpha}}\}=0 $, $|Z_{i,n,\alpha}| \leq 2$, and
$$
b_{n}^2 \eqdef \sum_{i=1}^n \Em \{ Z_{i,n,\alpha}^2 \}
\leq
n \Em \{ Q_{1,n,\alpha}^2 \} \leq 2 n F_{\nu}(\gamma_n).
$$
Since the $Z_{i,n,\alpha}$ are independent, we can apply
Bernstein's inequality.
If $\tau_{n} \eqdef \frac{\von}{2} n F_{\nu}(\gamma_n)$,
equation (6.9) and Bernstein's inequality entail:
$$
\Pm \bigl\{ \D_{n,\gamma_n,\alpha,K,\von}^c \bigr\}
\leq
2\exp \biggl( \frac{-\tau_{n}^2}{ 2(b_{n}^2 + 2 \tau_n /3) }\biggr)
\leq
2 \exp \Bigl( -\frac{\von^2}{8(2 + \von /3)} n F_{\nu}(\gamma_n) \Bigr),
$$
thus (6.7) follows. The proof of equation (6.8) is similar.
When $\beta=-1$, we have $\nu(t) =  t^{-1} \ell_{\nu}(t)$.
Define $Z_{i,n} \eqdef Q_{i,n,0} -  \Efm\{Q_{i,n,0}\}$.
In view of equation (A.5) we have
$$
\lim_{n \raro +\infty} \frac{1}{F_{\nu}(\gamma_n)} \Em\{ Q_{i,n,0}\}
= \lim_{n \raro +\infty} \frac{2}{F_{\nu}(\gamma_n)}
\int_0^1 K(t/h) \ell_{\nu}(t)\, dt/t
= 2 K(0) > 0.
$$
Then for $n$ large enough one has
$$
\Bigl\{ \Bigl| \frac{K_{n,\gamma_n,0}}{n F_{\nu}(\gamma_n)}
- 2K(0) \Bigr| > \von \Bigr\}
\subset
\Bigl\{ \Bigl| \frac{1}{n F_{\nu}(\gamma_n)}
\sum_{i=1}^n Z_{i,n} \Bigr| > \von/2 \Bigr\}.
$$
The $Z_{i,n}$ are independent and centered and $|Z_{i,n}| \leq 2$.
Moreover, in view of Assumption~K we have as before
$b_n^2 \eqdef \sum_{i=1}^n \Em\{ Z_{i,n}^2 \}
\leq 2 n  F_{\nu}(\gamma_n)$ and using again the Bernstein
inequality we get~(6.8).
\qed
\medskip

{\bf Lemma 6.}
{\it Let Assumption~{\rm K} hold. Assume that
$\omega \in \RV(s)$ with $s > 0$, $\mu \in \mc R(x_0,\beta)$
with $\beta > -1$, and $\lba_{\beta, K}$ is defined by equation
{\rm(4.1)}.
We have $\lba_{\beta,K} > 0$ and for any
$0 < \von \leq \frac{1}{2}$ we can find an event
$\mc  A_{n,\von}$ such that for $n$ large enough
$$
\mc A_{n,\von} \subset \{ | \lba(\mc X_{h_n}^K) - \lba_{\beta,K} |
\leq \von \}
\cap
\{ | \lba(\mc X_{H_n}^K) - \lba_{\beta,K} | \leq \von \}
\cap
\Big\{ \Big| \frac{H_n}{h_n} - 1 \Big| \leq \von \Big\}
\tag 6.10
$$
and
$$
\Pm\{ \mc A_{n,\von}^c \}
\leq 4(k+2) \exp \big( - c_{\beta, \sigma, \von} r_n^{-2} \big),
\tag 6.11
$$
where $c_{\beta, \sigma, \von} > 0$.}
\medskip

{\it Proof}.
Since $\lba_{\beta,K}$ is the smallest eigenvalue of
$\mc X_{\beta}^K$, we have $\lambda_{\beta,K} > 0$, otherwise
defining $\mb p(y) = (1, y, \ldots, y^k)$
and since $\mc X_{\beta}^K$ is symmetric, we should have
$$
0 = \lambda_{\beta, K}
= \inf_{\norm{x}=1} \prodsca{x}{\mc X_{\beta}^K x}
= \prodsca{x_0}{\mc X_{\beta}^K x_0}
= \int_{-1}^1 \bigl(\trans x_0 \mb p(y) \bigr)^2 y^{\beta} K(y)\,dy,
$$
where $x_0 \neq 0$ is the normalized eigenvector associated to
the eigenvalue $\lba_{\beta,K}$ and where we used the fact that
$$
\lba(M) = \inf_{\norm{x} = 1} \prodsca{x}{M x},
\tag 6.12
$$
for any symmetric matrix $M$. Then $\forall y \in \supp K$ we
have $\trans x_0\mb p(y) = 0$, which leads to a contradiction
since $y \mapsto \trans x_0 \mb p(y)$ is a polynomial.
For any $h, \von > 0$ we introduce the events:
$$
\aligned
\A_{n, h, \von} &= \bigl\{|\lambda(\mc X_{h}^K) - \lambda_{\beta, K} |
\leq \von \bigr\},
\\
\B_{n, h,\alpha, \von} &= \Bigl\{ \Bigl| \overline{K}_{n,h,\alpha}
-  \frac{\beta+1}{2} K_{\alpha,\beta} \Bigr| \leq \von \Bigr\}.
\endaligned
\tag 6.13
$$
Using the characterization (6.12) we can easily prove that
$$
\bigcap_{\alpha=0}^{2k} \B_{n, h, \alpha, \von/(k+1)^2}
\subset \A_{n, h, \von}.
\tag 6.14
$$
Since
$$
\align
\overline{K}_{n,H_n,\alpha} - \overline{K}_{n,h_n,\alpha}
&=
\overline{K}_{n,H_n,\alpha} \Big( 1 - \frac{N_{n,H_n}}{N_{n,h_n}}
\Big( \frac{H_n}{h_n} \Big)^{\alpha} \Big)
\\
&\qquad
+ \frac{1}{N_{n,h_n}} \sum_{i=1}^n
\Big( \frac{X_i - x_0}{h_n} \Big)^{\alpha}
\Big( K \Big( \frac{X_i - x_0}{H_n} \Big)
- K \Big( \frac{X_i - x_0}{ h_n} \Big) \Big),
\endalign
$$
we have when $K$ is the rectangular kernel $K^R$,
$$
| \overline{K}_{n,H_n,\alpha} - \overline{K}_{n,h_n,\alpha} |
\leq
\Big| \frac{N_{n,H_n}}{N_{n,h_n}} \Big( \frac{H_n}{h_n}
\Big)^{\alpha} -1 \Big|
+ \frac{1}{2} \Big( \frac{H_n}{h_n} \vee 1 \Big)^{\alpha}
\Big| \frac{N_{n,H_n}}{N_{n,h_n}} -1 \Big|,
$$
and otherwise under Assumption~K
$$
| \overline{K}_{n,H_n,\alpha} - \overline{K}_{n,h_n,\alpha} |
\leq
\Big| \frac{N_{n,H_n}}{N_{n,h_n}}
\Big( \frac{H_n}{h_n}\Big)^{\alpha} -1 \Big|
+ \frac{N_{n,H_n}}{N_{n,h_n}} \Big(\frac{H_n}{h_n} \Big)^{\alpha}
\rho \Big| \frac{H_n}{h_n} -1\Big|^{\kpa}
+ \rho \Big| \frac{h_n}{H_n} -1 \Big|^{\kpa}.
$$
Let us introduce for $\von > 0$ the event
$$
\F_{n,\von} \eqdef \Big\{ \Big| \frac{N_{n,H_n}}{N_{n,h_n}} - 1 \Big|
\leq \von \Big\}.
$$
Then for a good choice of $\von_1 \leq \von$ we have
$|\overline{K}_{n, H_n, \alpha} - \overline{K}_{n,h_n,\alpha} |
\leq
\frac{\von}{2(k+1)^2}$ on the event
$\C_{n,\von_1} \cap \F_{n,\von_1}$ and since $K \leq 1$,
we have $K_{\alpha, \beta} \leq \frac{2}{\beta+1}$
and noting that $\D_{n,h,0,K^R,\von_1}
= \bigl\{ \bigl| \frac{N_{n,h}}{2 n F_{\nu}(h)} - 1 \bigr|
\leq \von_1  \bigr\}$, we have for any $\alpha \in \setN$
$$
\D_{n, h, 0, K^R, \frac{\von}{3(k+1)^2 + \von}}
\cap \D_{n, h, \alpha, K, \frac{\von}{3(k+1)^2 + \von}}
\subset
\B_{n, h, \alpha, \frac{\von}{2(k+1)^2}}.
$$
Using (6.14) we get for $\eta \eqdef \frac{2 \von}{3(k+1)^2 + 2 \von}$:
$$
\D_{n,h_n, 0, K^R,\eta}
\cap
\bigcap_{\alpha=0}^{2k}
\D_{n, h_n, \alpha, K, \eta} \subset \A_{n, h_n, \von}.
\tag 6.15
$$
We take $0 < \von_2 \leq \von_1$ such that
$\frac{(1+\von_2)^{\beta+3}}{ 1 - \von_2 } \leq 1 + \von_1$
(for $\von_1$ small enough). Since $h \mapsto N_{n,h}$ is
increasing we have
$$
\C_{n,\von_2} \subset \{ N_{n,(1 - \von_2) h_n} \leq N_{n,H_n}
\leq N_{n,(1 + \von_2)h_n} \},
$$
and in view of Lemma~1 we can take $0 < \von_3 \leq \von_2$ such that
$$
\D_{n, (1 - \von_2) h_n, 0, K^R, \von_3}
\cap
\D_{n, (1 + \von_2) h_n, 0, K^R, \von_3}
\subset \C_{n,\von_2}.
$$
Using (A.1) with the slowly varying function
$\ell_F(h) \eqdef F_{\nu}(h) h^{-(\beta+1)}$,
we have for $n$ large enough that uniformly in
$y \in [\frac{1}{2}, \frac{3}{2}]$
$$
(1 - \von_1) \ell_F(h_n) \leq \ell_F(y h_n)
\leq (1 + \von_1)\ell_F(h_n),
\tag 6.16
$$
in particular, for $y = 1 - \von_1$ and $y = 1 + \von_1$ we get
by the definition of $\von_2$ and since
$\von_3 \leq \von_2 \leq \von_1$:
$$
\D_{n, (1 - \von_2) h_n, 0, K^R, \von_3}
\cap \D_{n, (1 + \von_2) h_n, 0, K^R, \von_3}
\cap \D_{n, h_n, 0, K^R, \von_3}
\subset \F_{n, \von_1}.
$$
Then we define for
$\von_4 \eqdef \von_3 \wedge \frac{\von}{3(k+1)^2 + \von}$
the event
$$
\mc A_{n,\von} \eqdef \D_{n, (1 - \von_2) h_n, 0, K^R, \von_4}
\cap \D_{n, (1 + \von_2)  h_n, 0, K^R, \von_4}
\cap \D_{n,h_n, 0, K^R, \von_4}
\cap \bigcap_{\alpha=0}^{2k} \D_{n, h_n, \alpha, K, \von_4},
$$
which satisfies (6.10) in view of the previous embeddings.
Using inequality (6.7) in Lemma~5
and since $\von _4 \leq \von_2 \leq \von_1 \leq \frac{1}{2}$,
we get
$$
\Pm\{ \mc A_{n,\von}^c \} \leq 4(k+2)
\exp \Big( -\frac{2^{-(\beta+ 3)} \von_4 \sigma^2 }
{ 8 ( 2 + \von_4 / 3)} r_n^{-2} \Big),
$$
where we used (6.16) and (2.5).
\qed
\medskip

\qquad{\tensmc6.4. Proof of the lower bounds}
\smallskip

{\bf Lemma 7.}
{\it If there are two elements $f_0$ and $f_1$ of a class
$\Sigma$ such that the Kullback--Leibler distance between the
corresponding probabilities $\Prob_0$ and $\Prob_1$ satisfies
$\mc K(\Prob_0, \Prob_1) < Q < +\infty$ with
$|f_0(x_0) - f_1(x_0)| \geq 2 c r_n$ for some constant $c > 0$,
then the pointwise minimax risk $\mc R_{n}(\Sigma, \mu)$
over the class $\Sigma$ defined by {\rm(2.1)}
in the model {\rm(1.1)} satisfies{\rm:}
$$
\mc R_{n}(\Sigma, \mu) \geq C(c, Q, p) r_n,
$$
where $ C(c, Q, p) \eqdef \frac{c}{2^{1/p}} \bigl(e^{-Q}
\vee\frac{1-\sqrt{Q/2}}{2} \bigr)^{1/p}$.}
\medskip

This result is classical. It can be found in Tsybakov~(2003) with a
proof based on a reduction scheme with two hypotheses and inequalities
between the Kullback--Leibler distance and other probability
distances.
\medskip

{\bf Proposition 6.}
{\it Let $h_n$ be defined by {\rm(2.5)}, let $(\alpha_n)$ be a
sequence of positive numbers going to $+\infty$ and
$r_n = \omega(h_n)$.
If $\Sigma = \Sigma_{h_n,\alpha_n}(x_0, \omega)$ is the class
given by Definition~{\rm2}, we have
$$
\liminf_n r_n^{-1} \mc R_{n}(\Sigma, \mu) \geq C_{s,p}.
\tag 6.17
$$}

{\it Proof}.
We use Lemma~7. All we have to do is to find two functions
$f_{0,n}$ and $f_{1,n}$ such that:
\roster
\item"(1)"
there is some $0 < Q < +\infty$ such that
$\mc K(\Prob_0^n,\mbb P_1^n) \leq Q$;
\item"(2)"
$f_{0,n}, f_{1,n} \in \Sigma_{h_n,\alpha_n}(x_0,\omega)$;
\item"(3)"
$| f_{0,n}(x_0) - f_{1,n}(x_0) | \geq 2 c r_n$
for some constant $c > 0$.
\endroster
We choose the two following hypotheses:
$$
f_{0,n}(x) = \omega(h_n)\ind{|x-x_0| \leq h_n},
\qquad
f_{1,n}(x) = \omega(|x-x_0|) \ind{|x-x_0| \leq h_n}.
$$

(1) Since the $\xi_i$ are centered Gaussian of variance
$\sigma^2$ and independent of $\mf X_n$, we have:
$$
\mc K(\Prob_0^n, \Prob_1^n \mid \mf X_n)
= \frac{1}{2 \sigma^2}\sum_{i=1}^{n}
\bigl( f_{0,n}(X_i) - f_{1,n}(X_i) \bigr)^2,
$$
then in view of (2.5)
$$
\mc K(\Prob_0^n, \Prob_1^n) =
\frac{n}{2\sigma^2} \norm{f_{0,n} - f_{1,n}}_{L^2(\mu)}^2
\leq \frac n{\sigma^2}  \omega^2(h_n) F_{\nu}(h_n)
= \frac 12.
$$

(2) For $h \in [0,h_n]$, taking $P$ as the constant polynomial
equal to $\omega(h_n)$, we have that the continuity modulus of
$f_{0,n}$ is~$0$, and taking $P = 0$ we obtain that the
continuity modulus of $f_{1,n}$ is bounded by $\omega(h)$.
Moreover, for~$n$ large enough, we clearly have
$f_{0,n}, f_{1,n} \in  \mc U(\alpha_n)$ since
$\alpha_n \raro + \infty$.

(3) If we take $c = 1/2$, we have
$|f_{1,n}(x_0) - f_{0,n}(x_0)| =  \omega(h_n) = 2 c r_n$.
\qed
\medskip

\qquad{\tensmc6.5. Computations of the examples}.
For a given design density, we compute the minimax convergence rate
$r_n$ by first giving an equivalent as $n \raro +\infty$ of the
smallest solution $h_n$~of
$$
\omega(h) = \frac{\sigma}{\sqrt{n F_{\nu}(h)}},
$$
and then an equivalent of $r_n = \omega(h_n)$.
\medskip

6.5.1. {\it Regularly varying design example}.
In the regularly varying design case we find the equivalent of
$h_n$ using the following proposition.
\medskip

{\bf Proposition 7.}
{\it  Let $\gamma >0$ and $\alpha \in \setR$.
If $G(h) = h^{\gamma} (\log(1/h))^{\alpha}$, then we have{\rm:}
$$
G^{\laro} (h) \sim \gamma^{\alpha / \gamma} h^{1/\gamma}
(\log(1/h))^{-\alpha / \gamma}
\qquad\text{as}\quad h \raro 0^+.
$$}
\medskip

{\it Proof}.
When $\alpha = 0$, the result is obvious, hence assume
$\alpha \in  \setR \setminus \{ 0 \} $.
We look for $h$ such that $h^{\gamma} (\log(1/h))^{\alpha} = x$,
when $x > 0$ is small. If $\alpha > 0$, we define
$t = \log(h^{\gamma / \alpha})$, so this equation becomes
$$
t \exp(t) = - \gamma x^{1/\alpha} / \alpha,
\tag 6.18
$$
where $t \leq 0$. The equation (6.18)
has two solutions for  $x$ small enough, but they cannot be
written in an explicit way.
Then let us consider the Lambert function $W$ defined as the
function satisfying $W(z) e^{W(z)} = z$ for any $z \in \mbb C$.
See, for instance, Corless \etal  (1996) about this function.
We are only interested here in its real branches.
This function has two branches $W_0$ and $W_{-1}$ in $\setR$.
We denote by $W_0$ the one such that $W_0(0)=0$ and $W_{-1}$
the one such that $\lim_{h \raro 0^-} W_{-1}(h) = -\infty$.
The two solutions of (6.18) are then
$t_0 = W_{-1} ( -\gamma x^{1/\alpha} / \alpha )$ and
$t_1 = W_0 ( -\gamma x^{1/\alpha} / \alpha )$ and
$h_0 \eqdef \exp \bigl( \alpha W_{-1}
( -\gamma x^{1/\alpha} / \alpha ) /\gamma \bigr)$
is the smallest solution.
By definition of $W$ we have for $-1/e < x < 0$
and $a \in \setR$:
$e^{a W_{-1}(x)} = (-x)^a (-W_{-1}(x))^{-a}$, and since
$W_{-1}$ satisfies $W_{-1}(-x) \sim \log(x)$
as $x \raro 0^+$, we have
$h_0 = (\gamma x^{1/\alpha} / \alpha)^{\alpha / \gamma}
( -W_{-1}(-\gamma x^{1/\alpha} / \alpha) )^{-\alpha / \gamma}
\sim  \gamma^{\alpha / \gamma} x^{1/\alpha}
( \log(1/x) )^{-\alpha / \gamma }$ as $x \raro 0^+$.

When $\alpha < 0$, we proceed similarly. We have $t \geq 0$ and
(6.18) has a single solution
$t = W_0 ( -\gamma x^{1/\alpha} / \alpha )$, thus
$h \eqdef \exp  -\alpha W_0( -\gamma x^{1/\alpha}/\alpha)/\gamma)$.
By the definition of $W_0$ we have $\forall x >0$ and $a \in \setR$:
$e^{a W_{0}(x)} = x^a W_{0}^{-a}(x)$, and since $W_{0}$ satisfies
$W_{0}(x) \sim \log(x)$ as $x \raro +\infty$,
we find again $h \sim \gamma^{\alpha / \gamma} x^{1/\alpha}
(\log(1/x))^{-\alpha / \gamma}$ as $x \raro 0^+$.
\qed
\medskip

For the second example of regularly varying design, using
Proposition~7, we find that an equivalent
to the sequence $h_n$ defined by (2.5) is
$$
(1 + 2s + \beta)^{ (\alpha + 2 \gamma) / (1 + 2s + \beta) }
\Big( \frac{ \sigma }{r} \Big)^{2 / (1 + 2s + \beta)}
(n (\log n)^{\alpha + 2 \gamma})^{-1 / (1 + 2s + \beta)},
$$
and since $\omega(h) = r h^s (\log(1/h))^{\gamma}$, we find that
an equivalent of $r_n$ (up to a constant depending on $s, \beta,
\gamma, \alpha$) is
$$
\sigma^{2s/(1+2s+\beta)} r^{(\beta+1)/(1+2s+\beta)}
(n (\log n)^{\alpha - \gamma(1+\beta)/s})^{-s/(1+2s+\beta)}.
$$
The computation for the third example ($\beta = -1$)
is similar to the second example, since
$F_{\nu}(h) = (\log(1 / h))^{1 - \alpha}$.
\medskip

6.5.2. {\it $\Gamma$-varying design example}.
For the $\Gamma$-varying design example
$\nu(h) = \exp(-1 /h^{\alpha})$, we first use the fact that when
$\nu \in \GV(\rho)$, we have $F_{\nu}(h) \sim \rho(h) \nu(h)$ as
$h \raro 0^+$ (see Appendix).
Recalling that $\rho(h) = \frac{h^{\alpha + 1}}{\alpha}$,
we solve
$$
h^{1+2s+\alpha} \exp(-1/h^{\alpha}) = y_n,
\tag 6.19
$$
where $y_n \eqdef \sigma^2 \alpha / (r^2 n)$.

Defining
$t \eqdef h^{-\alpha}$, equation (6.19) becomes
$t^{-(1+2s+\alpha)/\alpha} \exp(-t) = y_n$, which we rewrite as
$x \exp(x)= \alpha/(1+2s+\alpha) y_n^{-\alpha/(1+2s+\alpha)}$
for $x \eqdef\alpha/(1+2s+\alpha) t$.
Then we have
$x = W_0 \bigl( \alpha /(1+2s+\alpha)
y_n^{-\alpha/(1+2s+\alpha)}\bigr)$, where $W_0$ is defined in
the proof of Proposition~7.
Using the fact that $W_0(x) \sim \log(x)$ as $x \raro +\infty$,
we get
$x \sim \frac{\alpha}{1+2s+\alpha} \log n$ as $n \raro +\infty$,
thus
$h_n \sim (\log n)^{-1/\alpha}$ and the result holds since
$r_n \eqdef r h_n^s$.
\medskip

{\bf Appendix A. Some Facts on Regular and $\Gamma$-Variation}

We recall here some results about regularly and $\Gamma$-varying
functions. The results stated in this section can be found in
Bingham \etal (1989), Geluk and de Haan~(1987),
and Senata~(1976).
\medskip

\qquad{\tensmc A.1. Regular variation}.
Let $\ell$ be a slowly varying function throughout
the following. An important result is that the property
$$
\lim_{h \raro 0^+} \ell(y h) / \ell(h) = 1,
\tag A.1
$$
holds {\it uniformly\/} for $y$ in any compact set in
$(0,+\infty)$.
Now if $R_1 \in \RV(\alpha_1)$ and $R_2 \in \RV(\alpha_2)$,
one has
\roster
\item"(1)"
$R_1 \times R_2 \in \RV(\alpha_1 + \alpha_2)$,
\item"(2)"
$R_1 \circ R_2 \in \RV(\alpha_1 \times \alpha_2)$.
\endroster
If $R \in \RV(\gamma)$ for $\gamma \in \setR\setminus\{ 0 \}$,
then as $h\raro 0^+$ we have
$$
R(h) \raro   \cases
0       &\text{if}\quad \gamma > 0, \\
+\infty &\text{if}\quad \gamma <0.
\endcases
\tag A.2
$$
The asymptotic behaviour of integrals of regularly varying
functions, usually called Abelian theorems, plays a key role
in the proofs.

$\bullet$
If $\gamma > -1$ we have
$$
\int_{0}^{h} t^{\gamma} \ell(t)\, dt
\sim (1 + \gamma)^{-1} h^{1+\gamma} \ell(h)
\qquad\text{as}\quad h \raro 0^+,
\tag A.3
$$
and, in particular, $h \mapsto \int_0^h t^{\gamma} \ell(t)\, dt
\in \RV(\gamma + 1)$. This result is known as the Karamata theorem.

$\bullet$
When $\gamma = -1$ and if $\int_0^{\eta} \ell(t) \frac{dt}{t}
< +\infty$ for some $\eta>0$, then $h \mapsto \int_0^{h} \ell(t)
\frac{dt}{t} \in \RV(0)$ and we have
$$
\lim_{h \raro 0^+} \frac{1}{\ell(h)} \int_{0}^{h} \ell(t)
\frac{dt}{t} = +\infty.
$$

$\bullet$
If $R$ is some positive monotone function such that
$h \mapsto \int_0^h R(t)\, dt$ belongs to $\RV(\gamma)$ for some
$\gamma \geq 0$, then $R \in \RV(\gamma-1)$.

$\bullet$
If $K$ is a function such that
$\int_0^1 t^{-\delta} K(t)\, dt < +\infty$ for some $\delta >0$,
then
$$
\int_0^1 K(t) \ell(th)\, dt
\sim \ell(h) \int_0^1 K(t)\, dt
\qquad\text{as}\quad h \raro 0^+.
\tag A.4
$$
Moreover, when $\int_0^{\eta} \ell(t) dt/t < +\infty$
for some $\eta>0$, and $K$ is such that $\forall t \geq 0$,
$|K(t) - K(0)| \leq \rho |t|^{\kappa}$ for some $\rho >0$
and $\kappa >0$, one has
$$
\int_0^1 K(t/h) \ell(t) dt/t \sim K(0) \int_0^1 \ell(t) dt/t
\qquad\text{as}\quad h \raro 0^+.
\tag A.5
$$

If $R$ is defined and bounded on $[0, +\infty)$, one can define
the generalized inverse as
$$
R^{\laro}(y) = \inf\{ h > 0 \text{ such that } R(h) \geq y \}.
\tag A.6
$$
If $R \in \RV(\gamma)$ for some $\gamma > 0$, then there exists $R^{-}
\in \RV(1/\gamma)$ such that
$$
R ( R^{-} (h) ) \sim R^{-} ( R(h) ) \sim h
\qquad\text{as}\quad h \raro 0^+,
\tag A.7
$$
and $R^{-}$ is unique up to an asymptotic equivalence.
Moreover, one version of~$R^{-}$ is~$R^{\laro}$.

If $(\delta_n)_{n \geq 0}$ and $(\lba_n)_{n \geq 0}$ are
sequences of positive numbers such that $\delta_{n+1} \sim
\delta_n$ as $n \raro +\infty$,
$\lim_n \delta_n = 0$, and if there is a positive and
continuous function $\phi$ such that for any $y > 0$
$$
\lim_n \lambda_n R(y \delta_n) = \phi(y),
\tag A.8
$$
then $R$ varies regularly.
\medskip
\pagebreak

\qquad{\tensmc A.2. $\Gamma$-variation}.
We describe now the properties of $\Gamma$-varying functions and
$\Pi$-varying functions. The results are due to de Haan.
The references are the same as for regular variation. All the
following results can be found therein.

The first result states that if $\nu$ is a function such that
(2.6) holds for all $y \in \setR$, then (2.6) holds uniformly
on each compact set in~$\setR$. If~$\rho$ is such that (2.6)
holds, then
$$
\lim_{h \raro 0^+} \rho(h) / h = 0.
\tag A.9
$$
The auxiliary function $\rho$ in definition (2.6)
is unique up to within an asymptotic equivalence and can be
taken as $h \mapsto \int_0^h \nu(t) dt / \nu(h)$.

The class $\GV(\rho)$ is closed under integration.
If $\nu \in\GV(\rho)$, then
$F_{\nu}(h) = \int_0^h \nu(t)\, dt \in \GV(\rho)$
and we have
$$
F_{\nu}(h) \sim \rho(h) \nu(h)
\qquad\text{as}\quad h \raro 0^+.
$$

We have seen that
the class of regularly varying functions $\RV$ is closed
under the operation of functional inversion. In the
case of $\Gamma$-variation, the inversion maps the class $\GV$
in another class of functions, namely the de Haan class~$\PV$.
\medskip

{\bf Definition 5} ($\Pi$-Variation).
A function $\nu$ is in the de Haan class $\PV$ if there exists a
slowly varying function $\ell$ and a positive real number $c$
such that
$$
\forall y > 0, \quad
\lim_{h \raro 0^+} (\nu(yh) - \nu(h)) / \ell(y) = c \log(y).
\tag A.10
$$
The class of functions $\nu$ satisfying (A.10) is denoted
by $\PV(\ell)$.
\medskip

$\bullet$
If $\nu \in \GV(\rho)$, then
$\ell = \rho \circ \nu^{\laro}$ is slowly varying and
$\nu^{\laro} \in \PV(\ell)$.

$\bullet$
If $\nu \in \PV(\ell)$ for some $\ell \in \RV(0)$,
then $\nu^{\laro} \in \GV(\rho)$
with $\rho = \ell \circ \nu^{\laro}$.

In both senses the inverses and their auxiliary functions are
asymptotically unique. The following inclusion tells us that
$\Pi$-variation can be viewed as a refinement of slow variation.
Actually,
any $\Pi$-varying function is slowly varying:
for any $\ell\in \RV(0)$ we have
$$
\PV(\ell) \subset \RV(0).
\tag A.11
$$
\medskip

{\bf Acknowledgement.}
I wish to thank my adviser Marc Hoffmann for helpful suggestions
and encouragements.
\medskip
\medskip\medskip


\centerline{\bf References}
\medskip
{\eightpoint
\roster\widestnumber\item{[54]}
\item"[1]"
N.~H.~Bingham, C.~M.~Goldie, and T.~J.~Teugels (1989),
{\it Regular Variation},
Encyclopedia of Mathematics and its  Applications,
Cambridge University Press.
\item"[2]"
W.~S.~Cleveland (1979),
{\it Robust locally weighted regression and smoothing
scatterplots},
J\. Amer\. Statist\. Soc.,
{\bf 74}, 829--836.
\item"[3]"
R.~M.~Corless, G.~H.~Gonnet, D.~E.~G.~Hare, and D.~J.~Jeffrey (1996),
{\it On the Lambert $W$ Function},
Adv\. Comput\. Math., {\bf 5}, 329--359.
\item"[4]"
L.~de Haan (1970),
{\it On regular variation and its application to the weak
convergence of sample extremes}, PhD Thesis, University of
Amsterdam, Mathematical Centre Tract {\bf 32}.
\item"[5]"
J.~Fan and I.~Gijbels (1996),
{\it Local Polynomial Modelling and Its Applications},
Monographs on Statistics and Applied Probability,
Chapman \& Hall, London.
\item"[6]"
J.~L.~Geluk and L.~de Haan (1987),
{\it Regular Variations, Extensions and Tauberian Theorems},
CWI Tract.
\item"[7]"
A.~Goldenshluger and A.~Nemirovski (1997),
{\it On spatially adaptive estimation of  nonparametric
regression},
Math\. Methods Statist., {\bf 6}, 135--170.
\item"[8]"
E.~Guerre (1999), {\it Efficient random rates for nonparametric
regression under arbitrary designs},
Personal communication.
\item"[9]"
E.~Guerre (2000),
{\it Design adaptive nearest neighbor regression estimation},
J\. Multivariate Anal., {\bf 75}, 219--244.
\item"[10]"
P.~Hall, J.~S.~Marron, M.~H.~Neumann, and D.~M.~Tetterington (1997),
{\it Curve estimation when the design density is low},
Ann\. Statist., {\bf 25}, 756--770.
\item"[11]"
I.~A.~Ibragimov and R.~Z.~Hasminski (1981),
{\it Statistical Estimation: Asymptotic Theory},
Sprin\-ger, New York.
\item"[12]"
J.~Karamata (1930),
{\it Sur une mode de croissance r\'eguli\`ere des fonctions},
Mathematica (Cluj), {\bf 4} 38--53.
\item"[13]"
V.~Korostelev and A.~B.~Tsybakov (1993),
{\it Minimax Theory of Image Reconstruction},
Springer, New York.
\item"[14]"
A.~Nemirovski (2000),
{\it Topics in Non-Parametric  Statistics},
In: {\it Ecole d'\'et\'e de probabilit\'es de Saint-Flour XXVIII -- 1998},
Lecture Notes in Mathematics, no.~1738, Springer, New York.
\item"[15]"
S.~I.~Resnick (1987),
{\it Extreme Values, Regular Variation and Point Processes},
Applied Probability, Springer.
\item"[16]"
E.~Senata (1976),
{\it Regularly Varying Functions},
Lecture Notes in Mathematics, Springer.
\item"[17]"
V.~G.~Spokoiny (1998),
{\it Estimation of a function with discontinuities via local
polynomial fit with an adaptive window  choice},
Ann\. Statist., {\bf 26}, 1356--1378.
\item"[18]"
C.~J.~Stone (1977),
{\it Consistent nonparametric  regression},
Ann\. Statist., {\bf 5}, 595--645.
\item"[19]"
C.~J.~Stone (1980),
{\it Optimal rates of convergence for nonparametric estimators},
Ann\. Statist., {\bf 8}, 1348--1360.
\item"[20]"
A.~B.~Tsybakov (1986),
{\it Robust reconstruction of functions by the local
approximation},
Problems of Inform\. Transmission, {\bf 22}, 133--146.
\item"[21]"
A.~B.~Tsybakov (2003),
{\it Introduction \`a l'estimation non-param\'etrique},
Springer.
\endroster}
\medskip
\rightline{[{\sl Received October\/} 2004;
{\sl revised January\/} 2005]}
\end